\documentclass[reqno]{amsart}
\usepackage{epsfig,amsmath,amsfonts,latexsym}
\usepackage{amsthm}
\usepackage[abs]{overpic}
\usepackage[usenames,dvipsnames]{xcolor}
\usepackage{palatino}
\usepackage[hyperfootnotes=false]{hyperref}
\usepackage{caption}
\usepackage{dsfont}
\usepackage{comment}
\usepackage{accents}

\usepackage{amssymb}
\usepackage{braket}
\usepackage{adjustbox}

\usepackage[utf8]{inputenc} % Required for inputting international characters
\usepackage[T1]{fontenc} % Output font encoding for international characters

\usepackage{graphicx}
\usepackage{mathtools}
\usepackage{cite}
\usepackage{pifont}
\usepackage{tikz-cd}
\hypersetup{
  colorlinks,
  citecolor=Blue,
  linkcolor=Black,
  urlcolor=Green}
\usetikzlibrary{calc, 3d}

\definecolor{pakistangreen}{rgb}{0.0, 0.4, 0.0}
\definecolor{defnn}{RGB}{5, 66, 114} %Definitions
\definecolor{nicered}{HTML}{BA060F}

\theoremstyle{plain}

\newtheorem{theorem}{Theorem}

\newtheorem{dummy}{anything}[section]

\newtheorem{lemma}[dummy]{Lemma}

\newtheorem{proposition}[dummy]{Proposition}

\newtheorem{corollary}[dummy]{Corollary}

\theoremstyle{definition}
\newtheorem{definition}[dummy]{Definition}

\newtheorem{remark}[dummy]{Remark}

\newtheorem{assumption}{Assumption}
\theoremstyle{remark}
\newcommand{\R}{\mathbb{R}}

\newcommand{\norm}[1]{\left \lVert #1 \right \rVert}
\newcommand{\supp}{\textnormal{supp }}

\def\R{\mathbb{R}}

\renewcommand{\indent}{\hspace{1em}}

\title{A relative Poincaré--Birkhoff theorem}
\author{Agustin Moreno, Arthur Limoge}

\address[A.\ Moreno, A.\ Limoge]{Heidelberg University, DE}

\email{agustin.moreno2191@gmail.com, arthur.limoge@outlook.com}

\date{}

\begin{document}

\maketitle

%{\centering\footnotesize  \textit{To H.\ Poincar\'e, who taught us much;\\ To A. Floer, who followed suit;\\ To C.\ Viterbo, now on his 60th birthday, who took the cue; \\ and to all those who stand on the Shoulders of Giants.}\par}

\begin{abstract}
In \cite{MvK}, the first author and Otto van Koert proved a generalized version of the classical Poincaré-Birkhoff theorem, for Liouville domains of any dimension. In this article, we prove a relative version for Lagrangians with Legendrian boundary. This gives interior chords of arbitrary large length, provided the twist condition introduced in \cite{MvK} is satisfied. The motivation comes from finding spatial consecutive collision orbits of arbitrary large length in the spatial circular restricted three-body problem, which are relevant for gravitational assist in the context of orbital mechanics. This is an application of a local version of wrapped Floer homology, which we introduce as the open string analogue of local Floer homology for closed strings. 

%We further give an application of \cite[Theorem A]{MvK} proving, under suitable conditions, existence of infinitely many orbits for Reeb flows adapted to open books with symplectically trivial monodromy. This is in stark contrast with the symplectic topology underlying the three-body problem, where the presence of smoothly but \emph{not} symplectically trivial monodromy allows for the existence of finitely many orbits, as the famous Katok examples show \cite{K73}.   
\end{abstract}

\tableofcontents

\section{Introduction} 

In \cite{MvK}, following ideas of the proof of the Conley conjecture by Ginzburg \cite{G10}, the authors proved a generalized version of the classical Poincaré-Birkhoff theorem, for the case of Liouville domains of arbitrary dimension. The motivation for this result was finding spatial orbits in the well-known spatial circular restricted three-body problem (SCR3BP), concerning the dynamics of a negligible mass moving in $\mathbb{R}^3$ under the gravitational influence of two large bodies (the primaries). In this article, we prove a version of this result, but for the case of Lagrangians with Legendrian boundary, giving us insight into orbits of \textit{collision}.

\indent The motivation comes from the following geometric observation. In \cite{MvKb}, the authors proved, for low energies, the existence of global hypersurfaces of section for the SCR3BP. These give a reduction of the degrees of freedom, reducing the ambient continuous dynamics on a $5$-fold to the discrete dynamics of the Poincar\'e return map on a $4$-fold. This completely and non-perturbatively generalizes Poincar\'e's early idea of finding an annulus globally transverse to the planar three-body problem's dynamics – the earliest example of a hypersurface of section.

\smallskip

\indent We now know that in the CR3BP, these sections come in $\mathbb{S}^1$-families, and are the fibers of a fibration $\Sigma\backslash B \to \mathbb{S}^1$ (an open book) on the corresponding $5$-dimensional level set $\Sigma$, where the binding $B$ is the $3$-dimensional planar problem, an invariant subset. Intuitively speaking, the dynamics for the spatial problem $\Sigma$ "winds around" the planar problem $B$, while being transverse to all fibers. One usually denotes $M = \textbf{OB}(P, \phi)$ whenever $M$ admits such a fibration, where $P$ denotes the page (or fiber), and $\phi$ the monodromy, which allows us to glue the pages to each other. 

\smallskip

\indent For the SCR3BP, the topology of a page is either $\mathbb{D}^\star\mathbb{S}^2$ (for subcritical energy, and near each primary) or the boundary connected sum $\mathbb{D}^\star\mathbb{S}^2\natural \mathbb{D}^\star\mathbb{S}^2$ (for energy slightly above the first critical value). In the first case, if the page is chosen appropriately, the collision locus is precisely the cotangent fiber over the north pole in $\mathbb{D}^\star\mathbb{S}^2$, as we shall see in §\ref{CollisionsSection}; while in the second case we have two such cotangent fibers, corresponding to collisions with either primary. These are Lagrangians with Legendrian boundary, where the boundary corresponds to planar collisions. Intersections between a Lagrangian and its image under the return map correspond precisely to consecutive collision orbits. The fact that the dynamics continues after collision, allowing consecutive collisions, is an artifact of the way we shall regularize at singularities (collisions): forcing the small mass to bounce back whenever it collides with one of the primaries.

\indent Consecutive collision orbits can thus be recast as Hamiltonian chords between the corresponding Lagrangians, and are hence generators of wrapped Floer homology. Since the wrapped Floer homology of a cotangent fiber is well-known to be the homology of the free loop space over the base space, as we shall recall in §\ref{CollisionsSection}, it is here infinite-dimensional, so that one expects infinitely many consecutive collision orbits in the three-body problem. 

\medskip

While, of course, collision orbits are \emph{not} themselves of interest from a practical perspective (navigational engineers would very much prefer to avoid crashing a multi-million dollar spacecraft against the surface of Mars), these can normally be perturbed to actual orbits that pass close to the large celestial body, and use them e.g.\ for gravitational slingshots to reach another target, thus saving large amounts of precious fuel. In what follows, inspired by this heuristic idea, we shall prove a general statement which proves existence of infinitely many such chords (Theorem \ref{thm:longintchords}), assuming the twist condition introduced in \cite{MvK} holds (which generalizes the classical Poincaré--Birkhoff twist condition on the annulus).

\indent The proof of our main theorem will require the introduction of a \emph{local} version of wrapped Floer cohomology for Lagrangians in Liouville manifolds, as well as the existence of a local-to-global spectral sequence, built from local homologies, and converging to global wrapped Floer homology. This is the main technical novelty of this paper, and is carried out in Sections \ref{sec:locWFHdefn} and \ref{section:LocWFHproperties}. 

\medskip

\textbf{Hamiltonian twist maps.} Let us recall the twist condition from  \cite{MvK}. Let $(W,\omega=\mathrm{d}\lambda)$ be a $2n$-dimensional Liouville domain, and consider a Hamiltonian symplectomorphism $\tau$ of $W$. Let $(\partial W,\xi := \ker \alpha)$ be the contact manifold at the boundary where $\alpha=\lambda\vert_{\partial W}$, and $\mathcal{R}_\alpha$ is the Reeb vector field of $\alpha$. Recall that $\tau$ is Hamiltonian if $\tau=\phi_H^1$, where $\phi_H^t$ is the Hamiltonian flow of some $H_t$. We denote by $X_{H_t}$ its Hamiltonian vector field, defined by $i_{X_{H_t}}\omega = -\mathrm{d}H_t$. The Liouville vector field $V_\lambda$ is defined via $i_{V_\lambda}\omega=\lambda$. 

\begin{definition}(\textbf{Hamiltonian twist map})\label{def:twistmap}
We say that $\tau$ is a \emph{Hamiltonian twist map} (with respect to $\alpha$), if $\tau$ is generated by an at least $\mathcal{C}^2$ Hamiltonian $H:W \times \mathbb{R} \rightarrow \mathbb{R}$ which satisfies $X_{H_t}\vert_{\partial W}=h_t\mathcal{R}_\alpha$ for some \emph{positive} and smooth function $h:\partial W \times \mathbb{R} \rightarrow \mathbb{R}^+$.
\end{definition}

We will also need the notion of strong index-definiteness, as introduced in \cite{MvK}, in order to tell apart boundary chords from interior ones. We call a strict contact manifold $(Y,\xi=\ker \alpha)$ \emph{strongly index-definite} if the contact structure $(\xi,\mathrm{d}\alpha)$ admits a symplectic trivialization $\epsilon$ with the property that there are constants $c>0$ and $d\in \mathbb{R}$ such that for every Reeb chord $\gamma:[0,T]\rightarrow Y$ of length (or Reeb action) $T=\int_0^T \gamma^*\alpha$ we have
    $$
    \vert\mu_{RS}(\gamma;\epsilon)\vert\geq c T+d.
    $$

\bigskip

\textbf{Long interior chords.} Let us fix some basic conventions. Given a Hamiltonian symplectomorphism $\tau: W\rightarrow W$ of a symplectic manifold $W$, and a Lagrangian $L\subset W$, a pair $(x,m)$ where $x \in L$, $m \in \mathbb{N}\backslash\{0\}$, and $\tau^m(x) \in L$ is called a \textit{chord} of order $m$ of $\tau$. (It is the same thing as a length $1$ Hamiltonian chord of $H_t^{\# m}$, where $H_t$ generates $\tau$). We call \textit{minimal order} of $x$ the smallest such $m$. 

\indent We also call a chord periodic of period $k$ if $\tau^k(x) = x$, and define its \textit{minimal period} as the smallest such $k$ (note that it may be different from its minimal order). Note that a periodic chord is naught but a finite collection of \textit{sub-chords} $(x_0, m_0), \dots, (x_l, m_l) = (x_0, m_0)$, where each $(x_i, m_i)$ is obtained from the previous one by applying $\tau^{m_i -1}$. A periodic chord of period $1$ is one that is fixed by $\tau$.

\begin{remark}\label{NonGenericHamRk}
    For a generic Hamiltonian, end points of chords are never starting points of chords (Lemma 8.2 of \cite{AS10}). Hence, the minimal order is the same thing as the order, and periodic chords of period $>1$ simply do not exist. However, since we are interested in systems which are not necessarily generic (e.g.\ the CR3BP), we do not have the luxury of making such assumptions.
\end{remark}

\begin{theorem}[\textbf{Long interior chords}]\label{thm:longintchords}

Suppose that $\tau$ is an exact symplectomorphism of a connected Liouville domain $(W,\lambda)$. Let $\alpha:=\lambda\vert_{\partial W}$, and $L\subset (W,\lambda)$ be an exact, spin, Lagrangian with Legendrian boundary. Assume the following:
\begin{itemize}
\item\textbf{(Hamiltonian twist map)} $\tau$ is a Hamiltonian twist map (Definition \ref{def:twistmap}).
\item \textbf{(Periodic chords)} There are finitely many periodic chords.
\item\textbf{(Index-definiteness)} If $\dim W\geq 4$, then assume $c_1(W)\vert_{\pi_2(W)}=0$, and $(\partial W, \alpha)$ is strongly index-definite. 
\item\textbf{(Wrapped Floer homology)} $HW^*(L)$ is infinite dimensional.
\end{itemize}

\smallskip

Then $\tau$ admits infinitely many interior chords with respect to $L$, of arbitrary large order, and which are not sub-chords of any periodic chord.
\end{theorem}

Note that we make no assumption on the non-degeneracy of chords. 

\begin{remark}\label{rk:thmA} The following remarks are in order:

\begin{enumerate}
\item \textbf{(Periodic chords)} The second assumption of course doesn't hurt: if there already are infinitely many periodic chords, then it is still a positive result. As we point out in Remark \ref{NonGenericHamRk}, a generic Hamiltonian will actually have no periodic chords at all; but since the theorem is motivated by a physical application, we start with a given, possibly degenerate Hamiltonian, and hence cannot assume genericity.

    \item\textbf{(Wrapped Floer cohomology)} We impose the assumptions $c_1(W)\vert_{\pi_2(W)}=0$, as well as the assumptions on exactness and spin of the Lagrangian in order to have a well-defined and integer-graded cohomology. As for the last assumption: we do not know of any Lagrangian with $0 < \dim HW^*(L) < \infty$ (see \cite{SG} for a discussion on this folk question). In particular, we know that the collision Lagrangian from the SCR3BP satisfies all these hypotheses (Lemma \ref{CollLagrIsAdmissibleLemma} and Corollary \ref{CollisionOrbitsExistenceCorollary}).

    \item \textbf{(Index-definiteness)} This condition is motivated by Lemma \ref{IndexDefLemma} and Remark \ref{IndexDefRemark}. It was shown in \cite[Lemma D.1]{MvK} $\mathbb{R}^4$ that it holds whenever $W$ is a convex domain in $\mathbb{R}^4$. In particular, for $W = \mathbb{D}^\star \mathbb{S}^2$, the pages of the SCR3BP open book.
    %\item\textbf{(Spatial circular restricted three-body problem)} As in \cite{MvK}, we have not been able to check the twist condition for the return map in the SCR3BP, and hence the above theorem does not (as of yet) give the existence of long spatial consecutive collision orbits for this problem.
    
    \item \textbf{(Rotating Kepler problem)} The rotating Kepler problem, i.e.\ the two-body problem in a rotating frame, is the integrable limit case we obtain from the circular restricted 3-body problem when we set one of the large masses to zero. We shall study the collision orbits for this problem in Appendix \ref{app:RKP}, which can be understood as chords of a Lagrangian $2$-disc. It turns out that there are infinitely many chords of every order $k\geq 1$, all of which are periodic if a suitable resonance condition is satisfied.
\end{enumerate}    
\end{remark}

\textbf{Acknowledgements.} For this work, the authors were supported by the Air Force Office of Scientific Research (AFOSR) under Award No.\ FA8655-24-1-7012, by the DFG under Germany's Excellence Strategy EXC 2181/1 - 390900948 (the Heidelberg STRUCTURES Excellence Cluster), and by the Sonderforschungsbereich TRR 191 Symplectic Structures in Geometry, Algebra and Dynamics, funded by the DFG (Projektnummer 281071066 – TRR 191). \

\medskip

\section{Preliminaries on wrapped Floer cohomology}\label{section:WFHprereqs}

In this section, we review the definition of wrapped Floer cohomology for Lagrangians with Legendrian boundary in Liouville domains. Our motivation is the introduction of \textit{local} wrapped Floer cohomology in the next section, as a version of local Floer cohomology for open strings.

\subsection{Setup.} Consider a \emph{Liouville domain} $(W,\lambda)$, i.e. $W$ is a compact, exact symplectic manifold with boundary, and the vector field $V$ defined by the equation $\iota_V \mathrm{d}\lambda =\lambda$, called the \emph{Liouville vector field}, is outward pointing along each boundary component of $W$. The $1$-form $\lambda$ is called the \emph{Liouville form}, and its restriction $\alpha$ to $\partial W$ is a contact form. The associated \emph{Reeb vector field} $\mathcal{R}_\alpha$ is defined via the equations $\mathrm{d}\alpha(\mathcal{R}_\alpha,\cdot)=0,$ $\alpha(\mathcal{R}_\alpha)=1$. Given a Liouville domain $(W,\lambda)$, its completion is 
$$
(\widehat W, \widehat \lambda):= (W,\lambda) \cup_\partial ([1,\infty) \times \partial W,r \alpha), 
$$
and we denote $\widehat \omega:=d\widehat \lambda.$ Such completions are often called Liouville manifolds. Writing $r$ the coordinate on $[1,+\infty)\times\partial W$, the Liouville vector field is given by $V = r\partial_r$. 

\medskip

\begin{definition}\label{AdmissibleLagrDefn}
    A Lagrangian $L$ in $W$ is called \emph{admissible} if:
\begin{itemize}
    \item $L$ is transverse to $\partial W$, and $\partial L=L\cap \partial W$ is a Legendrian, i.e.\ $\alpha\vert_{\partial L}=0$;
    \item $L$ is exact, i.e.\ $\lambda\vert_L = \mathrm{d}f$ for some $f : L \to \mathbb{R}$; and
    \item $L$ is cylindrical near $\partial L$, i.e.\ the Liouville vector field $V$ is tangent to $L$ near $\partial L$.
\end{itemize}
\end{definition}

From the last condition, we can define the completion of $L$ as
$$
\widehat L=L\cup_{\partial L}[1,\infty)\times \partial L,
$$

and state the following simple lemma:

\begin{lemma}\label{lem:Hhat} Let $L\subset (W,\lambda)$ be an exact Lagrangian with Legendrian boundary. Then, after attaching an arbitrary small collar neighbourhood to $W$ along $\partial W$, we can extend $L$ to an admissible Lagrangian $\widehat{L}$.
\end{lemma}

\begin{proof}
We just need to ensure that the extension of $L$ intersects the boundary transversely and that it is cylindrical near its boundary. For this, we attach a small collar neighbourhood $[0,\varepsilon]\times \partial W$ to $\partial W=\{0\}\times\partial W$ along which $\lambda=e^t\alpha$ for $\alpha=\lambda\vert_{\partial W}$, and consider the Lagrangian $L_0=[-\varepsilon,\varepsilon]\times\partial L$. We choose a Weinstein neighbourhood $\mathcal{N}\cong T^*L_0$ of $L_0$. Then $L\cap [-\varepsilon,0]\times \partial W$ will lie in this neighbourhood after perhaps shrinking $\varepsilon$, as the image of an exact $1$-form $\mathrm{d}H$. We then extend $H$ to $\widehat H$ along the collar so that the image of $\mathrm{d}\widehat H$ intersects $\mathcal{N}\cap \left(\{\varepsilon\}\times \partial W\right)$ transversely, obtaining the desired extension of $L$. 
% \indent To show that $\widehat{L}$ is still exact (at least at infinity), observe that since $\widehat{L}\backslash \text{int}(L)$ is cylindrical above the Legendrian $\partial L$, $\hat{\lambda}$ is necessary zero over it. So the primitive of $\hat{\lambda}$ along $\widehat L \backslash \text{int}(L)$ should be constant. We already had a primitive $f$ for $\lambda$ on the interior of $W$. We can simply deform it to a constant, using a cut-off function, and therefore we will have $\hat{\lambda}|_{\widehat{L}} = \mathrm{d}f$ outside of some neighbourhood $(1, 1+\delta)\times\partial W$ above the boundary, where we can make $\delta$ arbitrarily small. 
\end{proof}

\bigskip

Now, let us introduce our dynamical objects. Consider a time-dependent, periodic Hamiltonian $H:\widehat W\times S^1 \to \R$, which we assume to be at least $\mathcal{C}^2$, and define its Hamiltonian vector field $X_H$ via
$$
\iota_{X_H} \widehat\omega = -\mathrm{d}H.
$$
We denote by $\phi_H^t: \widehat W \rightarrow \widehat W$ the Hamiltonian flow of $X_H$, and by $\phi=\phi_H^1$ its time-$1$ map. 

\begin{definition}
    A \textit{Hamiltonian chord} of $\widehat{L}$ is an orbit $x : [0,1] \to \widehat{W}$ of $X_H$ satisfying $x(0), x(1) \in \widehat{L}$.
\end{definition}

These naturally correspond to points in $\widehat L \cap \phi(\widehat L)$ via the correspondence $x \mapsto x(1)$. A Hamiltonian chord is called \textit{degenerate} if the corresponding intersection between $\widehat{L}$ and $\phi(\widehat{L})$ is non-transverse; and it is called \textit{contractible} if the class $[x] \in \pi_1(\widehat W, \widehat L)$ is trivial. 

\medskip

For the sake of generality, we will define our cohomology theory using two Lagrangians $\widehat{L}_0$ and $\widehat{L}_1$. Then, Hamiltonian chords will be $X_H$-orbits $x$ satisfying $x(i) \in \widehat{L}_i$, and those will be in natural correspondence with points in $\phi(\widehat{L}_0)\cap\widehat{L}_1$. Let us now impose two more assumptions on our system.

\smallskip

While $H$ is allowed to be time-dependent, we will hide this dependency in the next few sections, to make notation lighter. Including it makes no difference in the proofs, except each chord now comes in an $\mathbb{S}^1$-family, and has a trivial degeneracy in the $\mathbb{S}^1$-direction, which we ignore. See §8 of \cite{KvK} for this formalism.
\smallskip

\begin{assumption}\label{NonDegAssumption}
    $H$ is non-degenerate (i.e.\ it has no degenerate chords). We only make this assumption for §\ref{section:WFHprereqs}. Our constructions in §\ref{sec:locWFHdefn} and §\ref{section:LocWFHproperties} will allow us to lift it. 
\end{assumption}

\begin{assumption}\label{QuadHamAssumption}
    Let $r$ denote the interval coordinate on the collar neighbourhood $[1,+\infty)\times\partial W$. Assume there exists $r_0 \geq 1$ such that, for any $r \geq r_0$, $H = h(r)$, for some $h : [r_0, +\infty) \to \mathbb{R}$.
\end{assumption}

The second assumption has the following useful corollary:

\begin{corollary}\label{HamVFOnCollarLemma}
    On $[r_0, +\infty)\times\partial W$, the Hamiltonian vector field is given by $X_H = h'(r) \mathcal{R}_\alpha$.

    \indent Furthermore, if we denote by $\Omega(\widehat{L}_0,\widehat{L}_1)$ the space of smooth paths between $\widehat{L}_0,\widehat{L}_1$, let $f_i:\widehat{L}\rightarrow \mathbb R$ satisfy $\lambda\vert_{\widehat{L}_i}=df_i$ with $f_i\equiv 0$ for $r\geq r_0$, and consider the wrapped action functional $$\mathcal{A}_H:\Omega(\widehat L_0,\widehat L_1)\rightarrow \mathbb R,$$$$\mathcal{A}_H(x)=f_1(x(1))-f_0(x(0))-\int_0^1 x^*\widehat{\lambda}+\int_0^1H(x(t))\mathrm{d}t,$$ Then for any chord $x$ lying in $
\{r\geq r_0\}$, we have $\mathcal{A}_H(x) = -rh'(r) + h(r)$.
\end{corollary}

\medskip

Now, there are two standard schemes to construct wrapped Floer cohomology; one using linear Hamiltonians, and one using quadratic Hamiltonians; which we both explain now, for they will both be relevant if we want to state our constructions in full generality.

\begin{enumerate}
    \item First, observe that if one picks $H$ to be linear at infinity (i.e.\ for $r \geq r_0, H = h(r) = ar + b$, for $a>0,b \in \mathbb{R}$), then $X_H = a \mathcal{R}_\alpha$, so that if we choose $a$ in such a way that it is not the length of any Reeb chord on $(\partial W, \alpha)$, then $H$ will have \textit{no Hamiltonian chords} on $[r_0, +\infty)\times\partial W$, so that all chords are contained in a compact region.

    \indent One can then prove a maximum principle, ensuring that Floer trajectories joining these chords will also be contained in a compact region, hence ensuring that the usual Lagrangian Floer cohomology $HL^{*}(\widehat{L}_0, \widehat{L}_1; J, H)$ is well-defined (where $J$ is any suitably chosen compatible and cylindrical almost complex structure). While this cohomology is independent of $J$, it has the disadvantage on being dependent on $H$. To remediate that, one takes a direct limit over all such Hamiltonians, and defines:

    \vspace{-0.3em}

    \begin{equation}\label{WFHLinearDefn}
        HW^{*}(\widehat{L}_0, \widehat{L}_1) := \varinjlim\limits_{\text{slope}(H) \to \infty} HL^{*}(\widehat{L}_0, \widehat{L}_1; J, H),
    \end{equation}
    which we call the wrapped Floer cohomology of $\widehat{L}_0$ and $\widehat{L}_1$ in $\widehat{W}$ (or of $L_0$ and $L_1$ in $W$, since the Liouville completion was artificial). See \cite{Rit12} or \cite{JtK} for details.

    \medskip

    \item Though the above definition allows to associate a Floer cohomology theory to Lagrangians in a Liouville manifold, it is sometimes lacking. For example, Hamiltonians $H$ coming from a physical problem (whose chords have an actual physical relevance) are often degenerate, and do not necessarily have the prescribed behaviour near the boundary. We can always extend $H$ to some $\widehat{H}$ on $\widehat{W}$, and we would like to get information about the Hamiltonian chords of $H$ via the Floer theory of $\widehat{H}$. As we have seen above, if $\widehat{H}$ is linear at infinity, this can often be achieved (see e.g.\ \cite{MvK}). If $\widehat{H}$ grows faster, the above construction does not work. Instead, we have:

    \end{enumerate}

    \begin{proposition}\label{QuadWFHProp}
        Let $\widehat H$ be a Hamiltonian on $\widehat W$ of one of the following forms:

        \smallskip

        \begin{itemize}
            \item $\widehat H$ grows quadratically at $\infty$, i.e.\ $\widehat H \sim h(r) = \frac{1}{2}ar^2$ for $r \geq r_0$, $a > 0$.

            \item $\widehat H$ grows super-quadratically at $\infty$, i.e.\ $h(r)/r^2 \to \infty$.
        \end{itemize}

        Then one can associate a wrapped Floer cohomology $HW^{*}(\widehat{L}_0, \widehat{L}_1; J, \widehat{H})$ to $\widehat H$, and the latter will be naturally isomorphic to $HW^{*}(\widehat{L}_0, \widehat{L}_1)$, as defined in (\ref{WFHLinearDefn}).
    \end{proposition}

        \underline{\textit{Proof.}} (This argument is inspired from §18 of \cite{Rit12}, where the same proof is done for symplectic homology, except we take a different approach to define $HW^{*}$ at the end.)

        \smallskip

        By our assumption on $\widehat{H}$, and by Corollary \ref{HamVFOnCollarLemma}, $\partial_r \mathcal{A}_H$ becomes strictly negative (and bounded away from $0$) for $r$ big enough, so that $\mathcal{A}_H$ becomes strictly decreasing, and eventually goes to $-\infty$, as $r \to \infty$. Now, choose $R \gg 1$, and consider the region $\widehat{W}_R := \widehat{W}\cap\{r \leq R\}$. Then:

        \vspace{-0.6em}

        \begin{equation*}
            \left\{\text{Hamiltonian chords of } \widehat H \text{ in } \widehat{W}_R\right\} = \left\{\text{Hamiltonian chords of } \widehat H \text{ with } \mathcal{A}_H \geq a\right\}
        \end{equation*}

        \smallskip

        for some $a \ll -1$. In particular, consider the following construction.

    \smallskip

            \begin{center}
                \begin{tikzpicture}
  % Axis
  \draw[->] (0,0) -- (5,0) node[right] {$r$};
  \draw[->] (0,0) -- (0,5) node[above] {};

  \draw[domain=0:2,smooth,variable=\x,blue] plot ({\x},{0.5*\x*\x}) node[right] {$H_R \equiv \widehat H$};
  \draw[domain=2:3.5,smooth,variable=\x,purple] plot ({\x},{2*\x - 2}) node[right] {$H_R = mr + b$};

  % Black point at (2,2)
  \fill (2,2) circle (2pt) node[above right] {};

  % Dotted line to x-axis from (2,2)
  \draw[dotted] (2,2) -- (2,0) node[below] {$R$};
\end{tikzpicture}
    \end{center}

We construct a new Hamiltonian from $\widehat H$ by chopping it off after $r \geq R$, and completing it linearly with slope $m = h'(R)$. This gives us a Hamiltonian $H_R$ which is linear at infinity. Hence, by (1), it has no Hamiltonian chords on $[R, +\infty)\times\partial W$, giving us:

\vspace{-0.8em}

\begin{equation*}
    \left\{\text{Hamiltonian chords of } \widehat H \text{ with } \mathcal{A}_H \geq a\right\} = \left\{\text{Hamiltonian chords of } H_R\right\}.
\end{equation*}

Hence, the chain complexes $CF_{\geq a}^{*}(\widehat{L}_0, \widehat{L_1}; J, \widehat H)$ and $CF^{*}(\widehat{L}_0, \widehat{L_1}; J, H_R)$ are not just isomorphic, but equal. And we have already seen that the cohomology of the second one is $HL^{*}(\widehat{L}_0, \widehat{L}_1; J, H_R)$. Hence, the cohomology of $CF_{\geq a}^{*}(\widehat{L}_0, \widehat{L_1}; J, \widehat H)$ is well-defined, and:

\vspace{-0.6em}

\begin{equation}
    HW_{\geq a}^{*}(\widehat{L}_0, \widehat{L}_1; J, \widehat{H}) \cong HL^{*}(\widehat{L}_0, \widehat{L}_1; J, H_R).
\end{equation}

Notice that decreasing the action bound $a$ corresponds to increasing $R$, and thus the slope of the Hamiltonians $H_R$. Hence, in the direct limit:

\vspace{-0.6em}

\begin{equation}\label{DirectLimitEq}
    \varinjlim\limits_{a \to -\infty} HW_{\geq a}^{*}(\widehat{L}_0, \widehat{L_1}; J, \widehat{H}) \cong \varinjlim\limits_{R \to \infty} HL^{*}(\widehat{L}_0, \widehat{L}_1; J, H_R) = HW^{*}(\widehat{L}_0, \widehat{L}_1)
\end{equation}

\smallskip

Now finally: notice that we have never argued that $HW^{*}(\widehat{L}_0, \widehat{L}_1; J, \widehat{H})$ actually is well-defined (and indeed, this is quite a technical argument, see §18 of \cite{Rit12}). However, we \textit{have} argued that each $HW_{\geq a}^{*}$ is well-defined, and so we can take the following as a \textit{definition}:

\vspace{-0.8em}

\begin{equation*}
    HW^{*}(\widehat{L}_0, \widehat{L}_1; J, \widehat{H}) := \varinjlim_a HW_{\geq a}^{*}(\widehat{L}_0, \widehat{L}_1; J, \widehat{H}) \cong HW^{*}(\widehat{L}_0, \widehat{L}_1).
\end{equation*}

\smallskip
 (If we already knew $HW^{*}(\widehat{L}_0,\widehat{L}_1;J,\widehat H)$ to be well-defined, then this would hold from standard homological algebra; justifying our definition. See \cite{Lim} for details.) \qed 

\bigskip

The way we presented them, the second definition of $HW^{*}$ seems but like a re-packaging of the first one. However, the fact that it allows us to assign wrapped Floer cohomology to a single Hamiltonian will be essential in the next sections; in particular, in §\ref{section:LocWFHproperties}, to construct a spectral sequence starting from the local cohomology of each chord, and converging to $HW^{*}(\widehat{L}_0, \widehat{L}_1$).

\section{Local wrapped Floer cohomology: Definition}\label{sec:locWFHdefn}

\textbf{Notation.} Like earlier, $(\widehat{W}, \hat{\lambda})$ is a Liouville manifold, and $\widehat{L}_0$ and $\widehat{L}_1$ are admissible Lagrangians. We choose a compatible and cylindrical $J$ (in particular, $J\partial_r = \mathcal{R}_\alpha$), and $H : \widehat{W} \to \mathbb{R}$ a $\mathcal{C}^2$ Hamiltonian, on which we enforce for now no assumption (in particular, it could be degenerate), except:

\begin{assumption}
    Hamiltonian chords of $H$ between $\widehat{L}_0$ and $\widehat{L}_1$ are \textit{isolated}.
\end{assumption}

\smallskip

In this section, we propose a definition of local wrapped Floer cohomology at a given chord $x$; which, like other local homology theories, will record dynamical information about $x$ (how it may bifurcate under deformation of $H$, for example). Together with Proposition \ref{prop:SS}, which shows the existence of a local-to-global spectral sequence, this effectively allows one to assign wrapped Floer cohomology to a \textit{degenerate} Hamiltonian on $\widehat W$, and still have it be isomorphic to $HW^{*}(\widehat{L}_0, \widehat{L}_1)$.

\begin{definition}[\textbf{Local Wrapped Floer Cohomology}]\label{LocWFHDefn}
    Let $x$ be a Hamiltonian chord of $(J,H)$ in $\widehat W$, which may be degenerate. Then, one defines $HW_\text{loc}^{*}(x)$ as follows:

    \begin{itemize}
        \item choose a small neighbourhood $\mathcal{U}$ of $x$ in $\widehat{W}$ whose closure $\bar{\mathcal{U}}$ intersects no other chords of $(J,H)$. Such a $\mathcal{U}$ is commonly called an \textit{isolating neighbourhood} of $x$.

        \item call a pair $(\tilde{J},\tilde{H})$ a \textit{perturbation} of $(J,H)$ if it belongs to some neighbourhood of $(J,H)$, in the space of Floer data, endowed with one's desired topology. Since a generic such perturbation will be non-degenerate (see §8 of \cite{AS10}), its Lagrangian Floer cohomology is locally well-defined, allowing us to set:

        \vspace{-0.8em}

        \begin{equation}\label{HWlocDefExpr}
            HW_\text{loc}^{*}(\widehat{L}_0, \widehat{L}_1, x) := HL^{*}(\mathcal{U}, \widehat{L}_0\cap\mathcal{U}, \widehat{L}_1\cap \mathcal{U}; \tilde{J},\tilde{H})
        \end{equation}

        \smallskip

        \hspace{-1.3em} We call this the \textit{local wrapped Floer cohomology} of the chord $x$ (and often simply denote it $HW_\text{loc}^{*}(x)$, when the choice of Lagrangians is non-ambiguous).
    \end{itemize}
\end{definition}

To show that this construction is well-defined, there are three things to check:

\begin{enumerate}
    \item first, we need to show that if $(\tilde{J},\tilde{H})$ is close enough to $(J,H)$, then for any two chords $y,z$ of $(\tilde{J},\tilde{H})$ in $\mathcal{U}$, we can ensure that the Floer trajectories joining them \textit{stay in} $\mathcal{U}$. Once that is done, we still need to show that $HL^{*}(\mathcal{U})$ is well-defined; which is not automatic since it is usually constructed for closed manifolds, and $\mathcal{U}$ is non-compact.

    \item then, we show that the right-hand side of (\ref{HWlocDefExpr}) is independent of the choice of $(\tilde{J},\tilde{H})$.

    \item finally, we show that it does not depend on the choice of isolating neighbourhood $\mathcal{U}$ either.
\end{enumerate}

\smallskip

The rest of this section will be dedicated to proving these three steps. 

\bigskip

\textit{\textbf{Step 1 (Locality).}} The proof of \textbf{(1)} will rely on two lemmas. One, by now relatively standard, lemma from \cite{CFHW}, which we will introduce second; and one lemma which we call the «Energy Separation Property», allowing us to separate local from non-local trajectories based on their energy. This will be a straightforward analytical estimate – some traces of which can already be found in §4 of \cite{Oh96}, in a relatively similar set-up. Despite its simplicity, we want to insist on its importance: both in the proof of Step 1, as we will see; but also later on in the construction of the local-to-global spectral sequence, in Proposition \ref{prop:SS}, where this property will be essential.

\begin{lemma}[Energy Separation Property]\label{EnergySepLemma}
    Let $x$ be a (possibly degenerate) Hamiltonian chord. For any $\varepsilon_1 > 0$, there exist an isolating neighbourhood $\mathcal{U}$ of $x$, as well as a $\mathcal{C}^1$-neighbourhood $\mathfrak{U}$ of $(J,H)$, and $\varepsilon_2  > 0$ such that for any $(\tilde{J},\tilde{H}) \in \mathfrak{U}$, and for any Floer trajectory $u$ of $(\tilde{J},\tilde{H})$ intersecting $\mathcal{U}$, we have:

    \smallskip

    \begin{minipage}{0.75\linewidth}
        \begin{center}
            \begin{enumerate}
            \item $\text{im}(u) \subset \mathcal{U} \implies E(u) < \varepsilon_1$

            \item $E(u) < \varepsilon_2 \implies \text{im}(u) \subset \mathcal{U}$
        \end{enumerate}
        \end{center}
        
    \bigskip
    \bigskip

    where 

    \vspace{-0.6em}

    \begin{equation*}
        E(u) := \displaystyle\int_{\mathbb{R}\times[0,1]} \bigg|\dfrac{\partial u}{\partial s}\bigg|^2\mathrm{d}s\wedge\mathrm{d}t
    \end{equation*}

    \smallskip
    
    is the \textbf{energy} of $u$.

    \end{minipage}\hfill
    \begin{minipage}{0.25\linewidth}
    \centering
        \begin{tikzpicture}
    % Draw the arrow pointing upwards
    \draw[-latex, line width=0.8pt] (0,0) -- (0,4) node[above] {Energy};
    
    % Draw the blue graduation
    \draw[defnn, line width=0.8pt] (-0.1,1) -- (0.1,1) node[right, black] {$\textcolor{defnn}{\varepsilon_2}$};
    
    % Draw the red graduation
    \draw[nicered, line width=0.8pt] (-0.1,3) -- (0.1,3) node[right, black] {\textcolor{nicered}{$\varepsilon_1$}};
\end{tikzpicture}
    \end{minipage}
\end{lemma}

(Note that if such an $\varepsilon_2$ exists, then we have $\varepsilon_2 \leq \varepsilon_1$ without loss of generality. Indeed, if we had $\varepsilon_2 > \varepsilon_1$, then there could exist no Floer trajectories with energy in $[\varepsilon_1, \varepsilon_2)$, for their existence would contradict \textbf{(1)} and \textbf{(2)}; so that we may as well set $\varepsilon_2 = \varepsilon_1$).

\clearpage

    \underline{\textit{Proof.}} \textbf{(1)}  

    \begin{itemize}
        \item First, we look for a $\mathcal{C}^1$ estimate for the the action of perturbed chords. More precisely, given $(\tilde{J},\tilde{H})$ a $\mathcal{C}^1$ perturbation of $(J,H)$, and $\tilde{x}$ a chord of $(\tilde{J},\tilde{H})$, one can show that:
        
        \vspace{-0.6em}

            \begin{equation}\label{ActionEstimate0}
            \big|\mathcal{A}_{\tilde{H}}(\tilde{x}) - \mathcal{A}_H(x)\big| < C\norm{\tilde{x} - x}_{\mathcal{C}^1} + \big|\tilde{H} - H\big|_{\mathcal{C}^0},
        \end{equation} where $\mathcal{A}_H$ is the action functional of $H$. This estimate relies on standard Lipschitz-continuity arguments; and the full calculation can be found in \cite{Lim}. Note that $\norm{\tilde{x} - x}_{\mathcal{C}^1}$ is controlled by the diameter of $\mathcal{U}$, and the $\mathcal{C}^1$ distance between $(J,H)$ and $(\tilde{J},\tilde{H})$, since $\dot{x} = X_H(x)$. So we can rewrite (\ref{ActionEstimate0}) as:

        \vspace{-0.6em}

        \begin{equation}\label{ActionEstimate}
            \big|\mathcal{A}_{\tilde{H}}(\tilde{x}) - \mathcal{A}_H(x)\big| < \tilde{C}\left(\text{diam}\,\mathcal{U} + \norm{(J,H) - (\tilde{J},\tilde{H})}_{\mathcal{C}^1}\right).
        \end{equation}

        \smallskip

        (\underline{\textit{Note:}} of course, since $\mathcal{U}$ needs to contain the chord $x$, it cannot for example be a very small ball. However, we understand $\mathcal{U}$ to be a slightly enlarged tubular neighbourhood of $x$; in which case, when we say «shrinking $\mathcal{U}$», we really mean shrinking it in the normal direction to $x$; and by $\text{diam}\,\mathcal{U}$, we mean the diameter in that direction).

        \item 
        
        Therefore, by shrinking both $\mathcal{U}$ and $\norm{(J,H) - (\tilde{J},\tilde{H})}_{\mathcal{C}^1}$, one can make the right-hand side of (\ref{ActionEstimate}) smaller than $\varepsilon_1/2$. In particular, given two chords $y$ and $z$ of $(\tilde{J},\tilde{H})$ in $\mathcal{U}$, we have:

        \vspace{-0.6em}

        \begin{equation*}
            |\mathcal{A}_{\tilde{H}}(y) - \mathcal{A}_{\tilde{H}}(z)| < \varepsilon_1/2 + \varepsilon_1/2 = \varepsilon_1
        \end{equation*}

        by the triangle inequality.

        \item Since $\mathcal{A}_{\tilde{H}}(y) - \mathcal{A}_{\tilde{H}}(z)$ gives the energy of a Floer trajectory from $y$ to $z$, and since any Floer trajectory in $\mathcal{U}$ must end in Hamiltonian chords (\textit{a fortiori}, also in $\mathcal{U}$), we have proven \textbf{(1)}.    
    \end{itemize}

    \medskip

    \textbf{(2)} This is an easy proof by contradiction, which is inspired from the one in \cite{Oh96}. Fix the isolating neighbourhood $\mathcal{U}$ obtained in \textbf{(1)}. Then, assume \textbf{(2)} does not hold, so that there exists a sequence $(J_n, H_n) \longrightarrow (J,H)$, as well as a sequence $(u_n)$ of $(J_n, H_n)$-Floer trajectories, such that:

    \vspace{-0.8em}

    \begin{equation*}
        E(u_n) < \varepsilon_n \longrightarrow 0, \text{ but } \exists \, p_n \in \mathbb{R}\times[0,1] \text{ s.t } u_n(p_n) \not\in \mathcal{U}.
    \end{equation*}

    By elliptic regularity, one can extract a converging subsequence and a limit $u_\infty$. Since $E(u_\infty) = 0$, $u_\infty$ is constantly equal to a chord $x_\infty$ of $(J,H)$ in $\mathcal{U}$. And by our contradictory assumption: 

    \vspace{-0.7em}
    
    \[
    \exists t_\infty \in [0,1], x_\infty(t_\infty)\not\in\mathcal{U}.
    \]

    \indent However, since we assume our neighbourhood $\mathcal{U}$ to be isolating, it only contains the chord $x$, so that $x_\infty \equiv x$ up to translation, which gives the desired contradiction. \qed 

\bigskip

\begin{lemma}[\cite{CFHW}]\label{CFHWlemma}
    Let $x$ be a potentially degenerate chord of $(J,H)$, and $\mathcal{U}$ an isolating neighbourhood. Then, for every open $\mathcal{V} \subset\joinrel\subset \mathcal{U}$, there exists a $\mathcal{C}^\ell$ neighbourhood $\mathfrak{U}$ of $(J,H)$ (for $\ell \geq 1$ any integer) such that, for any $(\tilde{J},\tilde{H})$ in $\mathfrak{U}$, we have:

    \begin{itemize}
        \item all chords of $(\tilde{J},\tilde{H})$ contained in $\mathcal{U}$ are already contained in $\mathcal{V}$.

        \item all Floer trajectories of $(\tilde{J},\tilde{H})$ contained in $\mathcal{U}$ are already contained in $\mathcal{V}$.
    \end{itemize}
\end{lemma}

    \underline{\textit{Proof.}} Both steps are proven by contradiction, in much the same way as \textbf{(2)} right above. \qed

\bigskip

With these two lemmas at hand, we are now ready to conclude Step 1. Indeed, recall that we want to prove, formally, that:

\begin{proposition}
    Let $x$ be a Hamiltonian chord. One can find a small enough isolating neighbourhood $\mathcal{V}$ of $x$ such that: given $(\tilde{J},\tilde{H})$ a sufficiently $\mathcal{C}^1$-close perturbation of $(J,H)$, and $y$ and $z$ chords of $(\tilde{J},\tilde{H})$ contained in $\mathcal{V}$, then Floer trajectories connecting $y$ and $z$ do \textnormal{not} exit $\mathcal{V}$.
\end{proposition}

\underline{\textit{Proof.}} Fix $\varepsilon_1 > 0$, and recall that the Energy Separation Property, Lemma \ref{EnergySepLemma}, gives an isolating neighbourhood $\mathcal{U}$ of $x$, and some $\varepsilon_2 > 0$ such that:

\vspace{-0.6em}

\begin{equation}\label{TempEq1}
    E(u) < \varepsilon_2 \implies \text{im}(u) \subset \mathcal{U}.
\end{equation}

\smallskip

By the action estimate (\ref{ActionEstimate}), one can find a smaller isolating neighbourhood $\mathcal{V} \subset\joinrel\subset \mathcal{U}$, and choose $(\tilde{J},\tilde{H})$ sufficiently close to $(J,H)$ so that

\vspace{-0.6em}

\begin{equation*}
    E(u) < \varepsilon_2,
\end{equation*}

\smallskip

for any Floer trajectory $u$ of $(\tilde{J},\tilde{H})$, as long as its end chords are contained in $\mathcal{V}$. By (\ref{TempEq1}), $\text{im}(u) \subset \mathcal{U}$ for any such trajectory. We can now apply Lemma \ref{CFHWlemma}, which tells us that actually, $\text{im}(u) \subset \mathcal{V}$ (provided that we have chosen $(\tilde{J},\tilde{H})$ sufficiently close to $(J,H))$. Hence, every Floer trajectory $u$ with ends in $\mathcal{V}$ has image in $\mathcal{V}$. \qed

\medskip

In summary, we have shown that around every Hamiltonian chord $x$ of $(J,H)$, one can find a small isolating neighbourhood $\mathcal{V}$ of $x$ such that, for every sufficiently $\mathcal{C}^1$-close perturbation $(\tilde{J},\tilde{H})$ of $(J,H)$, all chords in $\mathcal{V}$ are joined by trajectories in $\mathcal{V}$. This makes it possible to define:

\vspace{-0.6em}

\begin{equation}
    HL^{*}(\mathcal{V}, \widehat{L}_0, \widehat{L}_1; \tilde{J},\tilde{H})
\end{equation}

The usual Lagrangian Floer cohomology of $\widehat{L}_0$ and $\widehat{L}_1$ in $\mathcal{V}$.

\begin{remark}\label{OpennessRemark1}
    There is actually a slight subtlety in writing the above expression: $\mathcal{V}$ is open, whereas Lagrangian Floer theory is usually defined on closed manifolds. However, recall that we obtained $\mathcal{V}$ as a pre-compact subset of $\mathcal{U}$ ($\mathcal{V} \subset\joinrel\subset \mathcal{U}$), and that Lemma \ref{CFHWlemma} ensured that all the Floer information (chords $\&$ trajectories) in $\mathcal{U}$ was already contained in $\mathcal{V}$. In particular, all the Floer information in $\mathcal{U}$ is bounded away from its boundary. Therefore, we can define $HL^{*}(\mathcal{U}, \widehat{L}_0, \widehat{L}_1; \tilde{J}, \tilde{H})$ without running into compactness or boundary issues. And since $HL^{*}(\mathcal{U}, \widehat{L}_0, \widehat{L}_1; \tilde{J}, \tilde{H})$ ultimately records Floer information about $\mathcal{V}$, we understand it to also mean $HL^{*}(\mathcal{V},\widehat{L}_0,\widehat{L}_1; \tilde{J},\tilde{H})$, and may use both notations interchangeably.
\end{remark}

This concludes Step 1.

\bigskip

\textit{\textbf{Step 2 (Perturbation invariance).}} The goal of this step is to show that the cohomology defined in Step 1 is independent of the choice of perturbation $(\tilde{J},\tilde{H})$, as long as the latter is sufficiently $\mathcal{C}^1$-close to $(J,H)$. This proof is the standard one, common to any Floer theory, and which can for example be found in Chapter 11 of \cite{AD10}, with little to no modification (the only 'difference' will be explained in the remark below). Let us very briefly sketch the argument (and refer to \cite{Lim} for a more detailed explanation of how we adapt it):

\medskip

Let us restrict ourselves to a $\mathcal{C}^1$-small neighbourhood of perturbations $\mathfrak{U}$ of $(J,H)$, and pick $(J_1, H_1)$ and $(J_2, H_2)$ in $\mathfrak{U}$. The idea is to show that, given a homotopy $\Gamma$ joining $(J_1,H_1)$ and $(J_2,H_2)$, one can construct a \textit{continuation map}:

\vspace{-0.6em}

\begin{equation}
    f_\Gamma : CF^{*}(J_1,H_1) \longrightarrow CF^{*}(J_2,H_2),
\end{equation}

\smallskip

in a similar fashion as to how one defines the Floer differential $\mathrm{d}$. Then, the same arguments which told us that $\mathrm{d}^2 = 0$ now tell us that $f_\Gamma$ is a chain morphism, and thus descends to a map on cohomology. Then, one can standardly deduce from the properties of these continuation maps (Ch. 11 of \cite{AD10}) that, if we denote by $\Gamma^{-}$ the time-reversal of $\Gamma$ (i.e.\ $\Gamma^{-}(s) = \Gamma(-s)$), then $f_\Gamma$ and $f_{\Gamma'}$ are inverses on cohomology. In our scenario, this gives:

\vspace{-0.6em}

\begin{equation}
    HL^{*}(\mathcal{V}, \widehat{L}_0, \widehat{L}_1; J_1, H_1) \cong HL^{*}(\mathcal{V}, \widehat{L}_0, \widehat{L}_1; J_2,H_2).
\end{equation}

\begin{remark}
    As claimed, the proof in our set-up works in exactly the same way as for Hamiltonian Floer theory; the details of which are fully fleshed out in Chapter 11 of \cite{AD10}. The only difference is that, since we work over an open manifold $\mathcal{V}$, arguments involving compactness may fail (when it comes to finding upper bounds, or limits to sequences of chords/trajectories).

    \indent This is easily fixable though. Indeed, since $\mathcal{V}$ is pre-compact in a bigger isolating neighbourhood $\mathcal{U}$, it suffices to choose some compact set $\mathcal{V} \subset K \subset \mathcal{U}$. One can then use upper bounds from $K$ to bound objects in $\mathcal{V}$. As for convergence: a sequence of chords/trajectories in $\mathcal{V}$ will have a converging subsequence in $K$. Lemma \ref{CFHWlemma} then tells us that this limit is actually contained in $\mathcal{V}$.
\end{remark}

This concludes the proof of Step 2.

\medskip

\textit{\textbf{Step 3 (Invariance on $\mathcal{V})$.}} We have hence successfully defined $HW_\text{loc}^{*}(x, \mathcal{V})$, and showed that it did not depend on the choice of Floer data $(J,H)$. It remains to show that, if we choose an even smaller isolating neighbourhood $\mathcal{V}' \subset\joinrel\subset \mathcal{V}$, then the induced cohomologies $HL^{*}(\mathcal{V})$ and $HL^{*}(\mathcal{V}')$ are isomorphic. This directly follows from Step 2. Indeed, pick some compact set $\mathcal{W} \subset \mathcal{V}'$, and choose a smooth bump function $\psi : \mathcal{V} \to \mathbb{R}$ such that:

\vspace{-0.6em}

\begin{equation*}
        \begin{cases}
            \psi|_\mathcal{W} \equiv 1 \\
            \psi|_{\mathcal{V}\backslash\mathcal{V}'} \equiv 0
        \end{cases}
\end{equation*}

\medskip

Now define $(J_\star, H_\star) := (\tilde{J},\tilde{H})\circ\psi$. Then $(\tilde{J}|_\mathcal{V},\tilde{H}|_\mathcal{V})$ and $(\tilde{J}|_{\mathcal{V}'}, \tilde{H}|_{\mathcal{V}'})$ can both be homotoped to $(J_\star,H_\star)$; and one can make sure that the second homotopy stays in $\mathcal{V}'$. By Step $2$, this gives us:

\vspace{-0.2em}

\begin{center}
     \begin{tikzcd}
	{HL^{*}(\mathcal{V}') \cong HL^{*}(\mathcal{V}'; \tilde{J}|_{\mathcal{V}'}, \tilde{H}|_{\mathcal{V}'})} & {HL^{*}(\mathcal{V}'; J_\star, H_\star)} & {HL^{*}(\mathcal{V}; \tilde{J}|_\mathcal{V}, \tilde{H}|_\mathcal{V})\cong HL^{*}(\mathcal{V})}
	\arrow["\cong", from=1-1, to=1-2]
	\arrow["\cong"', from=1-3, to=1-2]
\end{tikzcd}
\end{center}

\smallskip

which is the desired result. \qed

\bigskip

We have now shown Steps 1, 2, and 3, hence concluding that:

\vspace{-0.6em}

\begin{equation*}
    HW_\text{loc}^{*}(x) := HL^{*}(\mathcal{U}, \widehat{L}_0, \widehat{L}_1; \tilde{J},\tilde{H})
\end{equation*}

\smallskip

is well-defined, given any small enough isolating neighbourhood $\mathcal{U}$, and generic pair of non-degenerate Floer data $(\tilde{J},\tilde{H})$. 

\medskip

We have hence defined the local Wrapped Floer cohomology of any chord $x$.

\begin{remark}
    This is virtually a definition of local \textit{Lagrangian} Floer theory, since we are only working locally, and hence not detecting the global geometry of the Liouville manifold. However, we insist on calling it local \textit{wrapped} Floer cohomology because as we are about to see, in our set-up, one can recover \textit{global} Wrapped Floer homology from the local ones.
\end{remark}

\smallskip

\section{Local wrapped Floer cohomology: Properties}\label{section:LocWFHproperties}

Our notation is as earlier: we have two admissible Lagrangians $\widehat{L}_0,\widehat{L}_1$ in a Liouville manifold $(\widehat{W},\hat{\lambda})$, and a pair of potentially degenerate Floer data $(J,H)$, on which we assume:

\begin{assumption}\label{QuadHamAss2}
    $H$ is either quadratic, or grows faster than a quadratic at infinity (as in Prop. \ref{QuadWFHProp}); and Hamiltonian chords of $H$ are isolated. Denote by $\mathcal{P}(H)$ the set of such chords. In order to define the grading in $HW$, we need assume that $c_1(W)=0$ and the Lagrangians are spin (see after the proof for more details).
\end{assumption}

\begin{proposition}[Local-to-global spectral sequence]\label{prop:SS}
    Let $H$ be a (potentially degenerate) Hamiltonian, and $W,L_0,L_1$ satisfying Assumption \ref{QuadHamAss2}.  Then, there exists a half-plane spectral sequence $(E_n)$ whose first page is given by the local wrapped Floer cohomology of chords of $H$, and which converges to the global wrapped Floer cohomology of $(\widehat{L}_0, \widehat{L}_1)$ in $\widehat{W}$. In other words:

    \vspace{-0.6em}

    \begin{equation*}
        \displaystyle\bigoplus_{p \in \mathbb{N}} E_1^{p,*} = \bigoplus_{x \in \mathcal{P}(H)} HW_\text{loc}^{*}(x), \,\, E_1^{*,*} \implies HW^{*}(\widehat{L}_0, \widehat{L}_1).
    \end{equation*}

    \smallskip

    If, instead of Assumption \ref{QuadHamAss2}, $H$ is simply linear at infinity, then this spectral sequence also exists, but converges to $HL^{*}(\widehat{L}_0, \widehat{L}_1; J,H)$ instead.
\end{proposition}

One can be more precise about the form of the $E_1$ page. If one assumes that all Hamiltonian chords have distinct action, then each even column ($p \in 2\mathbb{N}$) corresponds to the local cohomology of one distinct chord, while odd columns are zero.

\indent If we do not assume that all chords have distinct action, then columns are separated by action, i.e.\ if two chords $x_1, x_2$ have the same action, then the corresponding column $p \in 2\mathbb{N}$ is given by:

\vspace{-0.6em}

\begin{equation*}
    E_1^{p,*} = HW_\text{loc}^{*}(x_1) \oplus HW_\text{loc}^{*}(x_2), * \in \mathbb{Z}.
\end{equation*}

\bigskip
\medskip

\textit{\underline{Proof of Prop. \ref{prop:SS}.}} By assumption, Hamiltonian chords are isolated, hence there are countably many of them. Arrange them in a sequence $\{x_k\}$, which is ordered by height $r$ for large $k$. Assume $H$ satisfies Assumption \ref{QuadHamAss2}. Then as we have already argued in Prop. \ref{QuadWFHProp}, $\mathcal{A}_H$ eventually becomes strictly decreasing, and goes to $-\infty$ as $r \to \infty$. The fact that it is \textit{strictly} decreasing tells us in particular that actions of chords do not accumulate. 

\indent Now fix $\varepsilon_1 > 0$. Recall that the Energy Separation Property (Lemma \ref{EnergySepLemma}) gives us, for every chord $x_k$, an isolating neighbourhood $\mathcal{U}_k$, and a constant $\varepsilon_2^k \in (0, \varepsilon_1]$ such that, given any perturbation $(\tilde{J},\tilde{H})$ of $(J,H)$ with Floer trajectory $u$ which intersects $\mathcal{U}_k$:

\vspace{-1.2em}

\begin{align*}
    \text{im}(u) \subset \mathcal{U}_k &\implies E(u) < \varepsilon_1 \\
    E(u) < \varepsilon_2^k &\implies \text{im}(u) \subset \mathcal{U}_k.
\end{align*}

\smallskip

This motivates us to construct $\varepsilon_2^k$ neighbourhoods around every $A_k = \mathcal{A}_H(x_k)$:

\begin{center}
\begin{tikzpicture}
    % Draw the horizontal axis
    \draw[->] (0,0) -- (10,0) node[anchor=north west] {Action $\mathcal{A}_H$};
    
    % Draw and label the points on the axis with different amplitudes for the brackets
    \draw[fill] (2,0) circle [radius=0.05];
    \node[below] at (2,0) {$A_{k-1}$};
    \draw[red] (2-0.4,0.3) arc[start angle=90,end angle=270,radius=0.3];
    \draw[red] (2+0.4,0.3) arc[start angle=90,end angle=-90,radius=0.3];
    \draw[<->,red] (2-0.4,0.5) -- (2+0.4,0.5) node[midway, above] {$\varepsilon_2^{k-1}$};
    
    \draw[fill] (5,0) circle [radius=0.05];
    \node[below] at (5,0) {$A_{k}$};
    \draw[red] (5-0.2,0.25) arc[start angle=90,end angle=270,radius=0.25];
    \draw[red] (5+0.2,0.25) arc[start angle=90,end angle=-90,radius=0.25];
    \draw[<->,red] (5-0.2,0.45) -- (5+0.2,0.45) node[midway, above] {$\varepsilon_2^k$};
    
    \draw[fill] (8,0) circle [radius=0.05];
    \node[below] at (8,0) {$A_{k+1}$};
    \draw[red] (8-0.6,0.35) arc[start angle=90,end angle=270,radius=0.35];
    \draw[red] (8+0.6,0.35) arc[start angle=90,end angle=-90,radius=0.35];
    \draw[<->,red] (8-0.6,0.55) -- (8+0.6,0.55) node[midway, above] {$\varepsilon_2^{k+1}$};
\end{tikzpicture}
\end{center}

\smallskip

Which we can choose to not overlap, since actions do not accumulate (we can always pick some smaller $\varepsilon_2^k$, without loss of generality).

\smallskip

Define a new sequence $(a_p)_{p \in \mathbb{N}}$ by taking, in order, the bounds of each of these intervals, i.e.\:

\begin{equation}\label{apSequence}
    \text{For } p = 2k: \begin{cases}
    a_p &:= A_k - \frac{1}{2}\varepsilon_2^k \\
    a_{p+1} &:= A_k + \frac{1}{2}\varepsilon_2^k.
\end{cases}
\end{equation}

\smallskip

Since $(J,H)$ may be degenerate, we fix a $\mathcal{C}^1$-close perturbation $(\tilde{J},\tilde{H})$ of $(J,H)$, and define $CF^{*}(J,H)$ to be the Floer co-chain complex of $(\widehat{W}, \widehat{L}_0, \widehat{L}_1; \tilde{J},\tilde{H})$. Then one can define the action filtration:

\vspace{-0.8em}

\begin{equation*}
    F_p CF^{*}(J,H) := \left\{ x \in CF^{*}(J,H) \, \mid \, \mathcal{A}_H(x) \geq a_p\right\}, p \in \mathbb{N},
\end{equation*}
where the grading is given by the CZ-index (as now generators are non-degenerate).
By standard arguments, a filtration on a co-chain complex induces a spectral sequence (see \cite{BottTu}), whose first page is given by:

\vspace{-0.8em}

\begin{equation}\label{E1pageExpr0}
    E_1^{p,*} = H^*(F_{p+1} CF^{*}(J,H)/F_p CF^{*}(J,H), \mathrm{d}^*),
\end{equation}

\medskip

where $\mathrm{d}^*$ is induced by the standard Floer differential. Observe that $F_{p+1} CF^{*}(J,H)/F_p CF^{*}(J,H)$ is generated by chords in $CF^{*}(J,H)$ with action in $(a_p, a_{p+1}]$, allowing us to explicitly compute $E_1$:

\begin{itemize}
    \item if $p$ is even, then $(a_p, a_{p+1}] = (A_k - \frac{1}{2}\varepsilon_2^k, A_k + \frac{1}{2}\varepsilon_2^k]$. Take two chords $y$ and $z$ of $(\tilde{J},\tilde{H})$ with actions in this interval. If $y$ and $z$ are both in the neighbourhood $\mathcal{U}_k$, then Floer trajectories connecting them will also be contained in $\mathcal{U}_k$, by the Energy Separation Property; so that $E_1^{p,*}$ at least contains the local wrapped Floer cohomology of the chord $x_k$.

    This is actually all of it. Indeed, if we assume that only one of $y$ or $z$ is contained in $\mathcal{U}$, then the Energy Separation Property still ensures that the Floer trajectory is fully contained in $\mathcal{U}$; so that both $y$ and $z$ were, in the end. And if neither is contained in $\mathcal{U}$, then any trajectory connecting them will have energy $\geq \varepsilon_2^k$, which contradicts our assumption that $y$ and $z$ have action in $(A_k-\frac{1}{2}\varepsilon_2^k, A_k +\frac{1}{2}\varepsilon_2^k]$. Therefore, we get:

    \vspace{-0.6em}

    \begin{equation}
        \text{for } p \text{ even:  } E_1^{p,*} = \displaystyle\bigoplus_{\mathcal{A}_H(x) = A_{p/2}} HW_\text{loc}^{*}(x).
    \end{equation}

    \item if $p$ is odd, then $(a_p, a_{p+1}] = (A_k + \frac{1}{2}\varepsilon_2^k, A_{k+1} - \frac{1}{2}\varepsilon_2^{k+1}]$. Hence, the cohomology (\ref{E1pageExpr0}) tries to connect chords with actions \textit{outside} of the $\varepsilon_2$ intervals. Since we have chosen those $\varepsilon_2$ neighbourhoods to not overlap, then by the action estimate (\ref{ActionEstimate}), there exist no such chords, provided that $(\tilde{J},\tilde{H})$ is sufficiently $\mathcal{C}^1$ close to $(J,H)$. Hence:

    \vspace{-0.6em}

    \begin{equation}
        \text{for } p \text{ odd: } E_1^{p,*} = 0.
    \end{equation}

  \end{itemize}

  \smallskip

    By standard arguments for filtration spectral sequences (see \cite{BottTu}), $(E_n)$ will now converge to the total cohomology of $F_p CF^{*}(\tilde{J},\tilde{H})$, which is by definition $HW^{*}(\widehat{L}_0, \widehat{L}_1)$. 

    \indent In more mathematical terms, with $A_k = \mathcal{A}_H(x_k) \,\, (k \in \mathbb{N})$, we have:

    \vspace{-0.5em}

    \begin{equation*}
        E_1^{p,q} = \begin{cases}
            \displaystyle\bigoplus_{\mathcal{A}_H(x) = A_{p/2}} HW_\text{loc}^{q}(x), \, &\text{ for } p \in 2\mathbb{N} \\
            \hspace{2em} 0 &\text{ for } p \in 2\mathbb{N} + 1
        \end{cases} \implies HW^{*}(\widehat{L}_0, \widehat{L}_1).
    \end{equation*} 

\smallskip

    Hence the name '\textit{local-to-global spectral sequence}'. \qed

\begin{remark}
    Convergence is actually not so straightforward given our definition of $(E_n)$, since this spectral sequence is unbounded. However, this can be fixed easily, by introducing a sequence of intermediary, bounded spectral sequences, with limit $(E_n)$:

    \indent Pick $a \in \mathbb{R}$, and define $(E_n^a)$, the filtration spectral sequence on the \textit{filtered} Floer complex $CF_{\geq a}^{*}(J,H)$ (consisting only of chords with action $\geq a$). Its $E_1$ page will be the same as we derived above, except only with chords of action $\geq a$. In particular, since $\mathcal{A}_H \to -\infty$, $(E_n^a)$ will be bounded for every $a$, and hence trivially converge, with limit $HW_{\geq a}^{*}(\widehat{W}, \widehat{L}_0, \widehat{L}_1)$.

    \indent Then, observe that $E_1 = \varinjlim_a E_1^a$, as $a \to -\infty$. And by exactness of inverse limits, we have:

    \vspace{-1em}

    \begin{align*}
        E_\infty = \varinjlim\limits_a E_\infty^a &\cong \varinjlim\limits_a HW_{\geq a}^{*}(\widehat{W}, \widehat{L}_0, \widehat{L}_1) \\
        &= HW^{*}(\widehat{W}, \widehat{L}_0, \widehat{L}_1)
    \end{align*} 

    This is what we wanted to prove. See \cite{Lim} for full details.

    \smallskip

    Also note that, in the exact same way as we constructed the sequence $(E_n^a)$, we can construct a spectral sequence for any Hamiltonian $H$ which is linear at $\infty$. Its first page will be given by the local cohomology of chords of $H$ (which are all contained in a compact region of $\widehat{W}$), and since the sequence is bounded, it will naturally converge to $HL^{*}(\widehat{L}_0, \widehat{L}_1; J, H)$. \qed 
\end{remark}

\medskip

\ding{226} The existence of this spectral sequence is the main technical result we shall need for the proof of Theorem \ref{thm:longintchords}. For now, let us introduce a last few definitions and simple lemmas that will be useful for this proof; as well as formally introduce the physical set-up, in §\ref{CollisionsSection}.

\smallskip

For Theorem \ref{thm:longintchords}, we shall assume that $c_1(W) = 0$ and that $L$ is spin. These ensure (see Remark 4.7 of \cite{Rit12}, or §9 of \cite{AS10}) that one can define a $\mathbb{Z}$-grading on $HW$. The latter is explicitly constructed in §3.6 of \cite{JtK}. For short: the assumption that $c_1(W) = 0$ ensures that we can trivialize $TW$ in a unique way, allowing us to define the usual Robbin-Salamon index $\mu_\text{RS} \in \frac{1}{2} \mathbb{Z}$ for symplectic paths, which we offset by $n/2$ (where $n := \dim W/2$) to get our grading on chords:

\vspace{-0.8em}

\[
\mu_\text{CZ} = \mu_\text{RS} - n/2 \in \mathbb{Z}
\]

\smallskip

called the Conley-Zehnder index. We can also define another index:

\begin{definition}\label{MeanIndexDefn}
    A map $f : G \to G'$ of Lie groups is said to be a \textit{quasimorphism} if for every $\psi, \phi \in G$, $|f(\psi\phi) - f(\psi) - f(\phi)| < C$, where $C$ is a constant that depends only on $G$ and $G'$. Let $\stackrel{\resizebox{10mm}{0.9mm}{$\ \sim$}}{Sp(2n)}$ be the universal cover of the symplectic group $Sp(2n)$. Then, there exists a unique quasimorphism $\Delta : \, \stackrel{\resizebox{10mm}{0.9mm}{$\ \sim$}}{Sp(2n)} \longrightarrow \mathbb{R}$, called the \textit{mean index}, and which satifies:

    \begin{enumerate}
        \item Let $\Phi$ be a symplectic path. For any close enough non-degenerate perturbation $\tilde{\Phi}$ of $\Phi$, we have $|\Delta(\Phi) - \mu_\text{CZ}(\tilde{\Phi})| \leq n$.

        \item Assume $\Phi$ is non-degenerate. Then $\lim\limits_{T \to \infty} \dfrac{\mu_\text{CZ}(\Phi|_{[0,T]})}{T} = \Delta(\Phi)$.

        \item Assume $\Phi = \gamma$ is a \textit{loop}, and write $\gamma^k$ its $k$-iteration. Then, $\Delta(\gamma^k) = k\Delta(\gamma)$.
    \end{enumerate}

    \medskip

    We refer to §3 of \cite{GG15} for a more thorough exposition, and interpretation of the mean index; the only fact we are interested in is the following:

    \begin{lemma}\label{SuppMeanIndexLemma}
        Let $H$ be some $\mathcal{C}^2$ Hamiltonian, and $x$ a Hamiltonian chord. Then:

        \vspace{-0.8em}

        \begin{equation*}
            \supp HW_\text{loc}(x) \subset [\Delta(x) - n, \Delta(x) + n]
        \end{equation*}

        \smallskip

        where $\supp HW_\text{loc}(x) := \{i \in \mathbb{Z} \, \mid \, HW_\text{loc}^{i}(x) \ne 0\}$.
    \end{lemma}

    \textit{\underline{Proof.}} This directly follows from property $(1)$ of Def. \ref{MeanIndexDefn}, since $\mu_\text{CZ}$ is the grading on $HW$. \qed 
\end{definition}

\medskip

Now that we have set all of this up, let us move on to a concrete application of this local wrapped Floer theory, to a problem in celestial mechanics.

\section{Collision orbits in the circular restricted three-body problem}\label{CollisionsSection}

The goal of this section is to turn the «search for collision orbits in the spatial circular restricted three-body problem» into a Floer theoretical problem, so that we can then use our constructions to prove an existence result in the next section. We will very briefly recall the set-up, and the Moser regularization process, for they are essential to defining the \textit{collision Lagrangian}.

\indent Three bodies move in $\mathbb{R}^3$: the Earth ($E$) and the Moon ($M$), which moves in circles about each other, so that we can pick a rotating frame and assume them to have fixed positions $q_E,q_M$; and a satellite $(S)$, of mass $m = 0$, whose motion is governed by the Hamiltonian:

\vspace{-0.8em}

\[
H(q,p) = \dfrac{1}{2}\norm{p}^2 - \dfrac{m_E}{\norm{q - q_E}} - \dfrac{m_M}{\norm{q - q_M}} + q_1 p_2 - q_2 p_1
\]

\smallskip

on $T^\star(\mathbb{R}^3\backslash\{q_E,q_M\})$. The first thing to do is to regularize at collisions. Indeed, assume the satellite collides with the Earth ($q \to q_E$). By conservation of energy ($H \equiv c$ along the motion), $p$ necessarily goes to $\infty$. In particular, the energy hypersurface $H^{-1}(c)$ is non-compact.

\indent To fix this, first define the symplectic switch map $(q,p) \mapsto (p, -q)$ on $T^\star \mathbb{R}^3 \cong \mathbb{R}^3 \oplus \mathbb{R}^3$, and compactify the bottom copy of $\mathbb{R}^3$ by adding the point at $p = \infty$. This gives us new coordinates $(\xi, \eta)$ on $T^\star \mathbb{S}^3$ (which formally, are obtained from $(p,-q)$ by applying the inverse stereographic projection). In these coordinates, our original Hamiltonian can be re-written as $Q(\xi, \eta) = \frac{1}{2}f(\eta, \xi)^2 \norm{\eta}^2$ (see e.g.\ §4 of \cite{MvKb} and references therein), which is but a deformation of the standard geodesic flow on $\mathbb{S}^3$. In particular, for low energies ($H \equiv c < H(L_1)$ where $L_1$ is the first Lagrange critical point), the bounded components of the regularized energy hypersurfaces of the SCR3BP (i.e.\ the level sets of $Q$) are diffeomorphic to $\mathbb{S}^\star\mathbb{S}^3$. This process is called \textit{Moser regularization}.

\indent By \cite{CJK18}, the above diffeomorphism is actually a contactomorphism, where $\mathbb{S}^\star\mathbb{S}^3$ is endowed with its standard contact structure, and the dynamics is given by a Reeb flow. In \cite{MvKb}, the authors explicitly construct an open book adapted to the dynamics, with pages $P \cong \mathbb{D}^\star\mathbb{S}^2$. This open book has the particularly nice feature that "its binding is the planar problem". 

\indent In the original phase space $T^\star(\mathbb{R}^3\backslash\{q_E,q_M\})$, the planar CR3BP (where the satellite is only allowed to move in the Earth-Moon plane) is naturally embedded as $$T^\star(\mathbb{R}^2\backslash\{q_E,q_M\}) \hookrightarrow T^\star(\mathbb{R}^3\backslash\{q_E,q_M\}).$$ This carries along through Moser regularization, giving an embedding $\mathbb{S}^\star\mathbb{S}^2 \hookrightarrow \mathbb{S}^\star\mathbb{S}^3$, where $\mathbb{S}^2$ sits inside $\mathbb{S}^3$ as the equator. Topologically, one then fixes the first page of the open book $P_0$, to be the $\mathbb{D}^\star\mathbb{S}^2$ obtained by those directions in the level set pointing into the upper-hemisphere of $\mathbb \mathbb{S}^3$ along this equator; and then gets the rest of the open book by pushing this page along the geodesic flow. We call the resulting open book the \emph{geodesic} open book. The actual open book found in \cite{MvKb} is more complicated, but it coincides with this one near the collision locus (defined below).

% \indent In particular, since the planar problem hypersurface $\mathbb{S}^\star\mathbb{S}^2 = \partial P_0$ is invariant under the geodesic flow (a planar trajectory does not suddenly become spatial), it remains the boundary of every page; thus the binding of the open book.

\smallskip

\ding{228} Now, assume the satellite collides with the Earth $(q \to q_E)$. As we saw at the start, we then have $p \to \infty$. And the point $p = \infty$ corresponds to the North pole $N$ of $\mathbb{S}^3$ (obtained by compactifying $\mathbb{R}^3$), so that the fiber of the level set, diffeomorphic to $\mathbb{S}_N^\star\mathbb{S}^3 \cong \mathbb{S}^2$, represents the \textit{collision locus}.

\indent Intuitively: recall that we had used the switch map $(q,p) \mapsto (p,-q)$ before compactifying, so that each vector in the fiber $\mathbb{S}_N^\star\mathbb{S}^3$ represents a direction $-q$ along which the satellite comes crashing into the Earth – it hence makes sense that this collision locus would be a $2$-sphere.

\begin{center}
    \begin{adjustbox}{scale=1.2}
        \begin{tikzpicture}[
  point/.style = {draw, circle, fill=red, inner sep=0.7pt},
]
\def\rad{2cm}
\coordinate (O) at (0,0); 
\coordinate (N) at (0,\rad); 

% Draw the bottom sphere
\filldraw[ball color=white] (O) circle [radius=\rad];
\draw[dashed] 
  (\rad,0) arc [start angle=0,end angle=180,x radius=\rad,y radius=5mm];
\draw
  (\rad,0) arc [start angle=0,end angle=-180,x radius=\rad,y radius=5mm];
\node at (0,0) {$\mathbb{S}^3$}; % Label for the bottom sphere

\begin{scope}[xslant=0.5,yshift=\rad,xshift=-2]
\filldraw[fill=gray!10,opacity=0.2]
  (-4,1) -- (3,1) -- (3,-1) -- (-4,-1) -- cycle;
\node at (2,0.6) {$T_N^\star\mathbb{S}^3$}; % Caption for the plane
\end{scope}
\node[point] at (N) {};
\node[above, red] at (N) {$N$};

% Draw a transparent smaller sphere centered at point N
\begin{scope}[shift={(0,\rad)}, opacity=0.3]
\def\smallrad{\rad/4}
\shade[ball color = blue!30] (0,0) circle (\smallrad);
\draw[dashed] 
  (\smallrad,0) arc [start angle=0,end angle=180,x radius=\smallrad,y radius=1.25mm];
\draw
  (\smallrad,0) arc [start angle=0,end angle=-180,x radius=\smallrad,y radius=1.25mm];
\node[above, blue] at (0, \smallrad+0.2) {$\mathbb{S}_N^\star\mathbb{S}^3$}; % Caption for the smaller sphere
\end{scope}

\end{tikzpicture}
    \end{adjustbox}
\end{center}

\smallskip

One can further observe that since $\mathbb{S}^\star\mathbb{S}^2 \hookrightarrow \mathbb{S}^\star\mathbb{S}^3$ is an energy hypersurface for the \textit{planar} three-body problem, the fiber over $N$ inside the level set, diffeomorphic to $\mathbb{S}_N^\star\mathbb{S}^2 \cong \mathbb{S}^1$, is the collision locus of the planar problem. It is embedded inside the spatial collision locus as the equator.

\begin{remark}[Adding transfer between the Earth and Moon]
    The scenario we have been describing thus far is for low energies $(c < H(L_1)$), in which case we have two such hypersurfaces $\mathbb{S}^\star\mathbb{S}^3$: one above the Earth, one above the Moon (in position space, they project down to the Hill regions). Taking the energy past the critical point $H(L_1)$ corresponds to adding a $1$-handle between these two hypersurfaces (a transfer window between the Earth and the Moon). Hence, for $c \in [H(L_1), H(L_1) + \varepsilon)$, we have one single hypersurface $\mathbb{S}^\star\mathbb{S}^3\natural \mathbb{S}^\star\mathbb{S}^3$. By \cite{CJK18}, it is still of contact-type (for small $\varepsilon$), and by \cite{MvKb}, there is still an open book decomposition, with pages $\mathbb{D}^\star\mathbb{S}^2\natural\mathbb{D}^\star\mathbb{S}^2$.

    \indent Hence, we now have two collision loci, which correspond to both fibers above the North poles (above the Earth and Moon). This will allow us to study situations where the satellite may collide with the Earth, bounce back, and collide with the Moon. In particular, such trajectories would tell us how to launch a satellite from Earth and send it crashing into the Moon. Admittedly, these are not trajectories we would like to use in practice. However, slight perturbations of them will give us near-miss trajectories, relevant for gravitational assist.
\end{remark}

\medskip

Now, recall that our objective is to use wrapped Floer theory in this set-up. And also recall that we have an open book on $\Sigma_c$ (the Moser regularized level set of energy $c$, which is the boundary of a fiberwise starshaped domain in $T^*\mathbb{S}^3$), whose pages are Liouville domains (up to the subtlety that the symplectic form degenerates at the boundary, which we will ignore). Therefore, let us look at the intersection of the collision locus with these pages. Denote by $N$ the North pole in $\mathbb{S}^3$. 

\begin{lemma}
    Let $P$ be any page of the open book in $\Sigma_c$. Then $P\cap T_N^\star \mathbb{S}^3$ is a Lagrangian submanifold of $P$.
\end{lemma}

\textit{\underline{Proof.}} Since $P \hookrightarrow \Sigma_c \hookrightarrow T^\star\mathbb{S}^3$, we can rewrite $P\cap \Sigma_c = P\cap T_N^\star\mathbb{S}^3$. Observe that $T_N^\star\mathbb{S}^3$ is a Lagrangian submanifold of $T^\star\mathbb{S}^3$. And since the open book is adapted to the dynamics, $\mathrm{d}\alpha|_P$ is symplectic along the interior of $P$ (where $\alpha$ is the contact form on the level set $\Sigma_c$). In particular, $\mathrm{d}\alpha|_{int(P)}$ is the restriction of the symplectic form $\omega = \mathrm{d}\lambda$ on $T^\star\mathbb{S}^3$ to $int(P)$, so that $P \cap T_N^\star\mathbb{S}^3$ is Lagrangian in $(int(P), \mathrm{d}\alpha|_{int(P)})$. \qed 

\medskip

We can refine this statement by choosing a particularly nice page of the open book. Pick $W := P_0$, the page of the ''geodesic'' open book in $\Sigma_c$ (constructed at the start of the section). Then, we have:

\smallskip

\begin{lemma}\label{CollLagrIsAdmissibleLemma}
    $L := W\cap T_N^\star\mathbb{S}^3$ is a Lagrangian disc in $W$, and it is admissible for wrapped Floer theory (Definition \ref{AdmissibleLagrDefn}). We call it the \textnormal{collision Lagrangian}.
\end{lemma}

\smallskip

\textit{\underline{Proof.}} Since the page $W = P_0$ agrees near the collision locus with the set of covectors along the equator pointing into the upper-hemisphere, then we directly have that $W\cap T_N^\star\mathbb{S}^3$ is identified with $\mathbb{D}_N^\star\mathbb{S}^2$, which is a disc, under the diffeomorphism $W\cong \mathbb D^*\mathbb \mathbb{S}^2$. Visually, one can view it as the upper-hemisphere of the collision locus $\mathbb{S}^2$ (the fiber of $\Sigma_c$ over $N$), whose boundary is the planar collision locus $\mathbb{S}^1$, embedded as the equator. 

\indent To show that $L$ is admissible for wrapped Floer cohomology (Definition \ref{AdmissibleLagrDefn}), we still need to show that $L$ is exact, and $\partial L := L \cap \partial W$ is Legendrian in $\partial W$. This is actually true of any fiber in a cotangent bundle. Hence, in particular, it holds for $\widehat{L} = T_N^\star\mathbb{S}^2$. \qed 

\medskip

Hence, $L$ is in the appropriate form to appeal to wrapped Floer cohomology. A Hamiltonian chord (i.e.\ a trajectory $x$ of the Hamiltonian vector field such that $x(0), x(1) \in L$) will correspond to an \textit{orbit of consecutive collisions} of the CR3BP. Indeed, recall that our Hamiltonian $H$ generates the Poincaré return map $\tau$ on the page $W$. Hence, a Hamiltonian chord of $L$ corresponds to a point $x \in L$ such that $\tau(x) \in L$. In other words: the satellite undergoes collision with the Earth, then bounces back, and collides with the Earth again. (\textit{Note}: this is possible because our change of coordinates to regularize collisions effectively made it so that, upon colliding with one of the primaries, the satellite bounced back where it came from).

\indent As in the Introduction, we call $x \in L$ a chord of order $m$ if $\tau^m(x) \in L$, and call $m$ its \textit{minimal order} if it is the smallest such $m$. And we call $x$ a \textit{periodic chord} if $\tau^k(x) = x$ for some $k$.

\begin{lemma}
    An interior chord of $L$ (contained in $\text{int}(W)$) corresponds to a spatial consecutive collision orbit, while a boundary chord of $L$ (contained in $\partial W$) corresponds to a planar consecutive collision orbit.
\end{lemma}

\textit{\underline{Proof.}} From the construction of our open book, the boundary of a page $W$ is the planar problem, so that a boundary chord is a Reeb chord with ends in the Legendrian $\partial L = L \cap \partial W$ (the planar collision locus, by Lemma \ref{CollLagrIsAdmissibleLemma}). Hence, such chords correspond to planar collision orbits. Meanwhile, an interior chord stays away from the binding, and is hence a spatial collision orbit. \qed 

\medskip

It remains to show that there is no chord with one end in the interior of $W$, and one end in the boundary. (This would correspond to a situation where the satellite undergoes spatial collision, bounces back, then somehow becomes tangent to the Earth-Moon plane to undergo a purely planar collision). But note that this scenario is forbidden since the planar problem is invariant.

\smallskip

We now make an assumption on the Hamiltonian generating the return map:

\begin{assumption}[Twist condition]\label{TwistConditionAssumption}
    The Poincaré return map $\tau : W \to W$ is a Hamiltonian twist map (Definition \ref{def:twistmap}), i.e.\ there exists a smooth positive (at least $\mathcal{C}^2$) function $h_t$ on $\partial W$ such that:

    \vspace{-0.8em}

    \[
    X_H|_{\partial W} = h_t \mathcal{R}_\alpha
    \]

    where $\mathcal{R}_\alpha$ is the Reeb vector field on $(\partial W, \mathcal{R}_\alpha)$, and $X_H$ is the Hamiltonian vector field of $H$, the Hamiltonian generating $\tau$.
\end{assumption}

\textit{\underline{Note:}} At this time, we cannot prove nor disprove that the return map for the SCR3BP satisfies the twist condition.

\medskip

This assumption has far-reaching consequences, which we will explore right after recalling the main result of this paper, which we have yet to prove:

\medskip

\textbf{Theorem \ref{thm:longintchords}.} \textit{Let $\tau$ be an exact symplectomorphism of a connected Liouville domain $(W, \lambda)$, and $L \subset (W,\lambda)$ be an exact, spin Lagrangian with Legendrian boundary. Further assume that:}

\begin{itemize}
    \item \textit{$\tau$ is a Hamiltonian twist map.}

    \item \textit{there are finitely many periodic chords}.

    \item \textit{if $\dim W \geq 4$, enforce that $c_1(W)|_{\pi_2(W)} = 0$, and ($\partial W, \alpha := \lambda|_{\partial W})$ is strongly index-definite.}

    \item $\dim HW^*(L) = \infty$.
\end{itemize}

\textit{Then, $\tau$ admits infinitely many interior chords on $L$, of arbitrarily large order, and which are not sub-chords of periodic chords.}

\begin{corollary}\label{CollisionOrbitsExistenceCorollary}
    Assuming the return map for the spatial circular restricted three-body problem satisfies the twist condition, then there exist infinitely many spatial orbits of consecutive collision, for the convexity range, and near the primaries.
\end{corollary}

Here, the convexity range is the set of all pairs of mass parameters and Jacobi constants for which the (Levi-Civita regularized) planar problem is strictly convex. By \cite{MvK}, strict convexity implies index-positivity for the dynamics.

\medskip

\bigskip

Before we move on with the proof of Theorem \ref{thm:longintchords}, let us motivate the second hypothesis, and the concept of «strong index-definiteness», introduced in \cite{MvK}. Recall from the Introduction:

\begin{definition}
    A strict contact manifold $(M, \xi := \ker \alpha)$ is said to be \textit{strongly index-definite} if there exists a symplectic trivialization $\epsilon$ of $\xi$, and constants $c > 0, d\in\mathbb{R}$ such that, for any Reeb arc $\gamma : [0,T] \to M$ (a portion of a Reeb trajectory):

    \vspace{-0.6em}

    \[
    |\mu_\text{RS}(\gamma;\epsilon)| \geq cT + d.
    \]
\end{definition}

\smallskip

The above definition can also be a made for a Hamiltonian flow on a symplectic manifold, where one now asks for the existence of a trivialization of the tangent bundle, satisfying the same condition. And we now have the following:

\begin{lemma}[Lemma 4.5 of \cite{MvK}]\label{IndexDefLemma}
    Let $(W, \lambda)$ be any Liouville domain, and $H : W \to \mathbb{R}$ satisfy the twist condition. Then, it can be extended to $\widehat{H} : \widehat W \to  \mathbb{R}$ such that:

    \begin{itemize}
        \item $\widehat H$ is linear on $[1+\delta,+\infty)\times\partial W$, where $\delta > 0$ can be made arbitrarily small.

        \item if $(\partial W, \alpha)$ is strongly index-definite, then so is the flow of $\widehat H$ on $\widehat{W}\backslash\text{int}(W)$.
    \end{itemize}
\end{lemma}

\textit{\underline{Proof sketch.}} First, the extension is obtained in a straightforward manner: one can write down the Taylor expansion of $H$ near $\partial W$ in the $r$-direction, and then artificially extend it to the whole collar $[1, +\infty)\times\partial W$. Then, one can use a cut-off function $\rho$ such that $\rho(r) \equiv 0$ for $r \geq 1 + \delta$, in order to cut off the terms of higher order, and turn $\widehat H$ linear. 

% More formally, we have:

% \vspace{-0.8em}

% \[
% H =  H_0 + (r-1)H_1 + R
% \]

% \smallskip

% where $H_0 := H|_{\partial W}, H_1 := (\partial_r H)|_{\partial W}$, and the remainder $R$ is $o(r-1)$. And we define:

% \vspace{-0.6em}

% \[
% \widehat{H} := \widehat{H}_0 + (r-1)\widehat{H}_1 + \widehat{R}
% \]

% where $\widehat{H}_i := \rho(r) H_i + (1-\rho(r))\max\limits_{\partial W} H_i$ and $\widehat{R} := \rho(r) R$.

% \medskip

\indent Hence, for $r \geq 1 + \delta$, we can write $\widehat H = ar + b$. In particular, $X_{\widehat{H}} = a \mathcal{R}_\alpha$, so that the flow of $\widehat H$ is trivially strongly index-definite on $[1+\delta, +\infty)\times\partial W$. The non-trivial part is then to show that this is also true on $[1, 1+\delta)\times\partial W$. This is the content of Lemma 4.5 of \cite{MvK}, where the authors linearize the Hamiltonian flow equation, to get $\dot{\psi} = \nabla X_{H} \, \psi$, and show that $\nabla X_H \in \mathfrak{sp}_{2n}$ is in the appropriate form to use a lemma from their Appendix A, stating that if the top-left $(2n-2)$-block of $\nabla X_H$ is already strongly index-definite, then so is the whole matrix. In particular, in our set-up, this top-left block will be asymptotically equivalent to $\nabla \mathcal{R}_\alpha$, the linearization of the Reeb vector field. This is already strongly index-definite by assumption, allowing them to conclude the proof. \qed 

\medskip

\begin{remark}[On Lemma \ref{IndexDefLemma}]\label{IndexDefRemark}
    The above lemma is of paramount importance for Theorem \ref{thm:longintchords}. Indeed, recall that we are interested in \textit{interior chords} of our Hamiltonian (spatial collision orbits). Since $\widehat H = ar + b$ is linear at infinity, we can safely assume that it has no chords on $[1+\delta,+\infty)\times\partial W$, by imposing that $a$ is not the length of any Reeb chord of $\partial L$ on the boundary.

    \indent However, we are still left with chords on the boundary $\partial W$ (which correspond to consecutive collision orbits in the \textit{planar} three-body problem, which we don't particularly care about in our set-up), and even worse: chords on $(1, 1+\delta)\times\partial W$, which are artifacts of our extension $H \to \widehat H$, and are therefore highly undesirable, for they have no physical significance.

    \indent The above lemma allows us to ignore such chords. Indeed, consider the sequence given by $\widehat H^{\# k}$, where $k \in \mathbb{N}\backslash\{0\}$, and $\widehat{H}^{\# k}$ is the $k$-th iterate of $\widehat H$. Then, a length $1$ Hamiltonian chord $x$ of $\widehat{H}^{\# k}$ corresponds to a length $k$ Hamiltonian trajectory of $\widehat{H}$ (note that there is no reason to assume that this will be a chord for $\widehat{H}$. In fact, generically, it won't be, by Remark \ref{NonGenericHamRk}; but we only need it to be a trajectory). As such, by strong index-definiteness:

    \vspace{-0.8em}

    \[
    |\mu_\text{RS}(x; \epsilon)| \geq ck + d.
    \]

    In particular, as we take $k \to \infty$, the Robbin-Salamon index of chords of $\widehat{H}^{\# k}$ on $\widehat{W}\backslash\text{int}(W)$ goes to infinity (and so does their Conley-Zehnder index, since $\mu_\text{RS}$ and $\mu_\text{CZ}$ only differ by a constant). Hence, these chords will not count towards the final cohomology $HW^{*}(L) = \varinjlim_k HW^{*}(L; \widehat{H}^{\# k})$.

    \indent Therefore, we can safely ignore chords on $\widehat{W}\backslash\text{int}(W)$, as we wanted.
\end{remark}

\smallskip

Let us state one last, simple lemma, before we move on to the proof of Theorem \ref{thm:longintchords}.

\begin{lemma}\label{lemma:HW_spread_out}
    Assume $HW^*(L)$ is infinite-dimensional, and $(\partial W, \alpha)$ is strongly index-definite. Then, $|\supp HW(L)| = \infty$, where $\supp HW(L) \subset \mathbb{Z}$ is defined as in Lemma \ref{SuppMeanIndexLemma}.
\end{lemma}

\textit{\underline{Proof.}} This is proven in Lemma 3.15 of \cite{MvK}, for symplectic homology. The argument is exactly the same here. \qed

\section{Proof of Theorem \ref{thm:longintchords}} In this section, we give a proof of Theorem \ref{thm:longintchords}, as an adaptation of the proof of \cite[Thm.\ A]{MvK}.

\begin{proof}[Proof of Thm.\ \ref{thm:longintchords}]
Write $\tau =\phi^1_H$, and let $L\subset W$ be a Lagrangian as in the statement of Thm.\ \ref{thm:longintchords}. We consider the extensions $\widehat W,\widehat H, \widehat{\tau}$ respectively of $W$, $H$, $\tau$ as in \cite{MvK}, so that $\widehat H$ is linear at infinity. We also consider the admissible extension $\widehat L$ to $\widehat W$ provided by Lemma \ref{lem:Hhat}.

\indent By assumption, $H$ has finitely many periodic chords, which are hence isolated. We denote by $c_1,\ldots,c_k$ the collection of fixed chords (i.e.\ the $x_j=c_j(0)$ are fixed points of $\tau$ lying in $L$). Also, we denote by $d_1,\dots,d_q$ the interior periodic chords with minimal period strictly greater than $1$, and by $n_1,\dots,n_q$ the corresponding minimal periods.

\indent Assume for a contradiction that $\tau$ has a finite number of interior chords which are not sub-chords of periodic chords. Denote their orders by $m_1 < \dots < m_\ell$. Take any sequence $\{p_i\}_{i=1}^\infty$ going to infinity, such that each $p_i$ is indivisible by $n_1,\dots,n_q$, and $p_i\geq \max_{j,k}\{n_j,m_k\}$ for every $i$. We may apply \cite[Lemma 4.6]{MvK} if necessary to perturb the Hamiltonian $\widehat H^{\# p_i}$ on the cylindrical part $[1,r_\infty]\times \partial W$, for some $r_\infty$. This ensures that $\widehat H^{\# p_i}$ is weakly admissible, so we can appeal to local wrapped Floer homology and the spectral sequence from Proposition~\ref{prop:SS} to compute $HL^*(\widehat L;\widehat H^{\# p_i})$.

\indent By Lemma~\ref{lemma:HW_spread_out}, for all $N>2n k$ we may find $i_1,\ldots,i_N$, ordered by increasing absolute value, and such that such that $HW^{i_j}(\widehat L)\neq 0$. Combining this with Lemma \ref{IndexDefLemma}, we may choose $p_i$ sufficiently large such that the following hold:
\begin{enumerate}
\item Each chord of $\widehat{H}^{\# p_i}$ that is contained in $\widehat W \setminus \mbox{int}(W)$ has index whose absolute value is larger than $|i_N|+2n$;
\item the Floer cohomology groups $HL^{i_j}(\widehat L,\widehat H^{\# p_i})$ are non-trivial for $j=1,\ldots, N$. 
\end{enumerate}
From the spectral sequence of Proposition~\ref{prop:SS} applied to $\widehat H^{\# p_i}$, and (2), we deduce that there are non-trivial summands on $E^1_{pq}(\widehat L,\widehat H^{\# p_i})$ with $p+q=i_j$ for $j=1,\ldots, N$. From (1) we know that no chord in $\widehat W \setminus \mbox{int}(W)$ can contribute to local Floer homology of degree $i_j$, and so we conclude that every term $E^1_{pq}(\widehat L,\widehat H^{\# p_i})$ in the spectral sequence with $p+q=i_j$ must come from the local wrapped Floer homology of a chord in int$(W)$.

\indent Because we have assumed that the $p_i$'s are indivisible by $n_1,\ldots, n_q$ we conclude that each such interior chord cannot be the iterate of one of the periodic chords $d_1,\ldots,d_q$. Moreover, since we chose the $p_i$ to be larger than $\max_k m_k$, we conclude that it must be an iterate $c_j^{p_i}$ of a fixed chord $c_j$. However, by (3) of Def. \ref{MeanIndexDefn}, we have $\Delta(c_j^{p_i}) = p_i\Delta(c_j)$, and by Lemma \ref{SuppMeanIndexLemma}:

\vspace{-0.8em}

\[
\mbox{supp}HW_{\text{loc}}^*(c_j^{p_i},\widehat H^{\# p_i})\subset[p_i \Delta(c_j)-n, p_i \Delta(c_j)+n].
\]

\smallskip

This covers at most $2n k$ different degrees, leaving some of the degree $i_j$ uncovered as we had chosen $N>2n k$. This is a contradiction. Therefore, there exist infinitely many interior chords which are not the sub-chord of any periodic chord.
\end{proof}

\appendix

\section{Collision orbits in the rotating Kepler problem}\label{app:RKP}

\begin{figure}
    \centering
    \includegraphics[width=0.9 \linewidth]{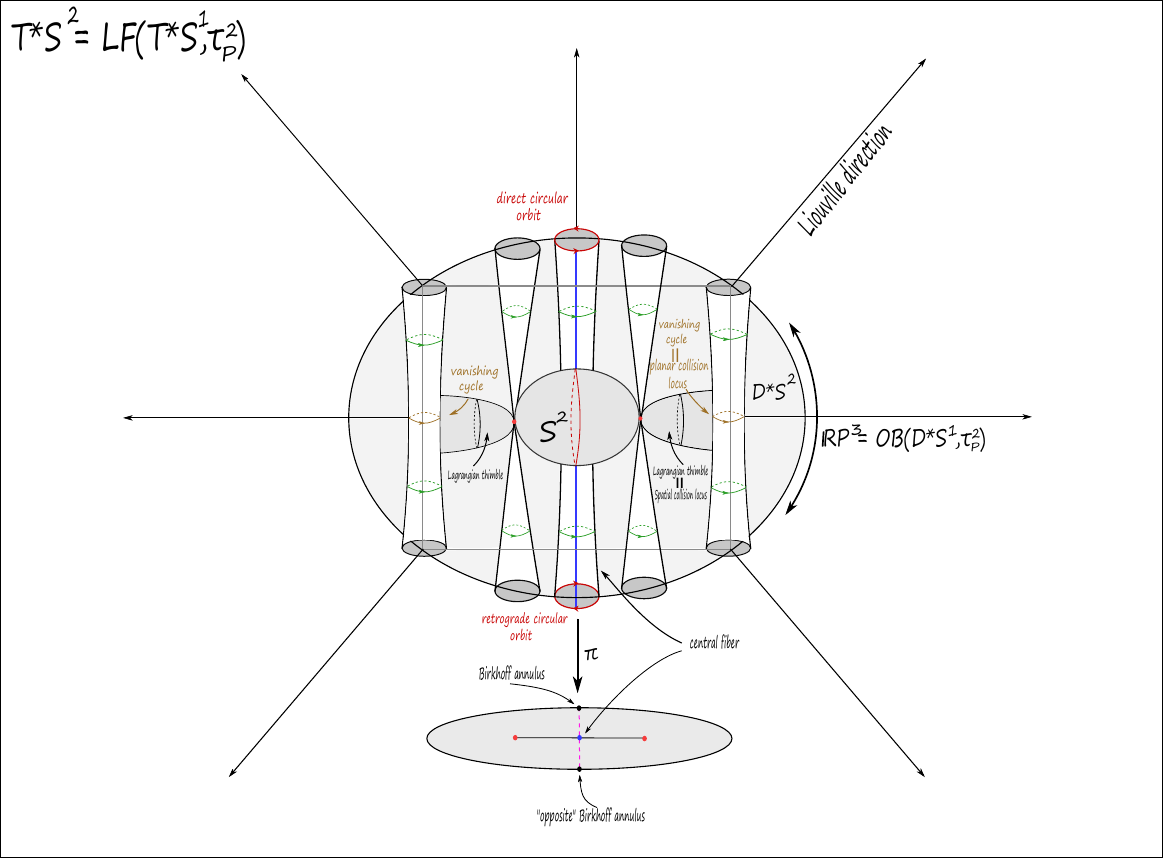}
    \caption{In the subcritical case of the SCR3BP, for suitably chosen page, the spatial collision locus corresponds to one of the Lagrangian thimbles of the standard Lefschetz fibration on $\mathbb{D}^*\mathbb{S}^2$ (shown above). Its boundary, the corresponding vanishing cycle, is the planar collision locus; this is a circle which can be identified with the zero section of a regular fiber. In the case where the planar problem admits an adapted open book (e.g.\ when the planar dynamics is dynamically convex \cite[Thm.\ 1.18]{HSW}), the existence of an adapted Lefschetz fibration as above was shown in \cite{Mo}. In this case, this vanishing cycle can be thought of as the zero section of an annuli-like global surface of section for the planar problem.}
    \label{fig:my_label}
\end{figure}

In this appendix, we study consecutive collisions for the rotating Kepler problem, which correspond to chords of the collision Lagrangian (a Lagrangian $2$-disc on the page $\mathbb{D}^*\mathbb{S}^2$ of a supporting open book). This is a completely integrable system, obtained as the limit of the circular restricted three-body problem by setting the mass of one of the primaries to zero. In unregularized coordinates, the dynamics is given by the Hamiltonian $H=K+L$, where
$$
K= \frac{1}{2} \Vert p \Vert^2  -\frac{1}{\Vert q\Vert},
\quad
L= q_1 p_2 - q_2 p_1,
$$
both of which are preserved quantities. If we write 
$$
T^*\mathbb{S}^3=\{(\xi=(\xi_0,\dots,\xi_3),\eta=(\eta_0,\dots,\eta_3))\in \mathbb{R}^4\oplus \mathbb R^4: \langle\xi, \eta\rangle=0, \Vert \xi \Vert=1\},
$$
then the Moser-regularized Hamiltonian is $Q:T^*\mathbb{S}^3\rightarrow \mathbb{R}$, $Q(\xi,\eta)=\frac{1}{2}f^2(\xi,\eta)\Vert \eta \Vert^2$ 
where
$f(\xi,\eta)=1+(1-\xi_0)(-c-1/2+\xi_2\eta_1-\xi_1\eta_2)$, with $c$ the Jacobi constant for $H$. 

\smallskip

It was shown in \cite[App.\ A]{MvKb} that the geodesic open book $$
(
\xi, \eta) \longmapsto \frac{\xi_3 +i\eta_3}{\Vert \xi_3+i\eta_3 \Vert}
$$
is a supporting open book for the rotating Kepler problem for $c<-3/2$. The Hamiltonian $L$ is given in these coordinates by $L=\xi_2 \eta_1 -\xi_1 \eta_2$. With respect to the $\pi/2$-page
$$
P= \left\{ (\xi;\eta) \in T^*\mathbb{S}^3:\;Q(\xi,\eta)=\frac{1}{2},\; \xi_3=0,\; \eta_3 \geq 0 \right\}\cong \mathbb{D}^*\mathbb{S}^2,
$$

\smallskip

the Poincaré return map is given by

\smallskip

$$
f: P\rightarrow P,
$$
$$
(\xi_0,\xi_1,\xi_2,0;\eta_0,\eta_1,\eta_2,\eta_3) \longmapsto
(\xi_0,R_{T(c-L)}(\xi_1,\xi_2),0;\eta_0,R_{T(c-L)}(\eta_1,\eta_2),\eta_3),
$$
where $R_\phi$ denotes the rotation by angle $\phi$, and $T(K)=\frac{\pi}{2(-K)^{3/2}}$ (the period of a Kepler ellipse of energy $K$). The collision locus is $\mathcal C=\{(\xi,\eta)\in P: \xi_0=1\}\cong \mathbb D^2$, i.e.\ the $2$-disc cotangent fiber over the north pole $N=(1,0,0,0)$ (also the north pole of $P\cong \mathbb{D}^*\mathbb{S}^2$). Using that $L\vert_{\mathcal C}=0$, we clearly see that $\mathcal C$ is invariant under $f$, and $f:\mathcal C\rightarrow \mathcal C$ is a rotation by angle $T(c)$, i.e.\ $f\vert_\mathcal{C}=R_{T(c)}$. 

We conclude that there are infinitely many chords of every order $k\geq 1$, although all of them have minimal order $1$. These are non-isolated, since they come in a family parametrized by the disc $\mathcal{C}$, whose boundary circle corresponds to planar chords. The origin is always fixed, and corresponds to the northern polar collision orbit. Moreover, if $T(c)/2\pi$ is irrational, there are no periodic chords except for the origin. If $T(c)=2\pi p/q$ is a rational multiple of $2\pi$, every point in $L$ different from the origin is periodic of the same minimal period $q$, and so gives a periodic chord of minimal period $q$.

\end{document}